%% file: Brinkman_paper.tex
\documentclass[review]{elsarticle}
\makeatletter
\def\ps@pprintTitle{%
 \let\@oddhead\@empty
 \let\@evenhead\@empty
 \def\@oddfoot{\centerline{\thepage}}%
 \let\@evenfoot\@oddfoot}
\makeatother

\pdfoutput=1
\usepackage{graphicx}
\usepackage{amsfonts, amsmath, amssymb}
\usepackage{booktabs,siunitx}
\usepackage[svgnames,table]{xcolor}
\usepackage[tableposition=above]{caption}
\usepackage{pifont}
\usepackage{bm}
\usepackage{cancel}
\usepackage{caption}
\usepackage{color}
\usepackage{enumerate}
\usepackage{float}
\usepackage{hyperref}
\usepackage[english]{babel}
\usepackage[margin=1in]{geometry}
\usepackage{mathtools}
\usepackage{multirow}
\usepackage{setspace}
\usepackage{subfigure}
\usepackage{lineno}
\newcommand \D [2]{\frac{\partial #1}{\partial #2}}

\renewcommand{\vec}[1]{\bm{\mathrm{#1}}}
\newcommand{\V}[1]{\bm{\mathrm{#1}}}

\def \div{\nabla \cdot \mbox{}}
\def \grad{\nabla}

\def \x{\vec{x}}

\def \qn{q^{\text{n}}}
\def \qnj{q^{\text{n}}_j}
\def \qd{q^{\text{d}}}
\def \qdi{q^{\text{d}}_i}
\def \n{\vec{n}}

\def \u{\vec{u}}

\def \e{\vec{e}}

\def \Nalign{N_{\rm aligned}}
\def \Nnalign{N_{\rm non-aligned}}

\def \vbeta{\boldsymbol{\beta}}

\def \fb{f_\text{b}}

\def \Omegas{\Omega_{\text{s}}}
\def \Omegaf{\Omega_{\text{f}}}

\def \e{\vec{e}}

\def \f{\vec{f}}

\def \half{\frac{1}{2}}
\def \3half{\frac{3}{2}}
\def \5half{\frac{5}{2}}

\def \n{\vec{n}}

\def \ncells{n_{\text{cells}}}

\def \pexact{p_\text{exact}}
\def \qexact{q_\text{exact}}

\def \u{\vec{u}}
\def \ub{\u_{\text{b}}}

\def \uexact{u_\text{exact}}

\def \vexact{v_\text{exact}}

\def \x{\vec{x}}

\def \div{\nabla \cdot \mbox{}}
\def \grad{\nabla}

\def \dt{\Delta t}

\def \dt{\Delta t}

\def \ndot{\n \cdot}

\newcommand{\upperRomannumeral}[1]{\uppercase\expandafter{\romannumeral#1}}

\begin{document}
\let\today\relax

\begin{frontmatter}
	
\title{Critique on ``Volume penalization for inhomogeneous Neumann boundary conditions modeling scalar flux in complicated geometry"}
\author[SDSU]{Ramakrishnan Thirumalaisamy}
\author[SDSU]{Nishant Nangia}
%\ead{nnangia@sdsu.edu}
\author[SDSU]{Amneet Pal Singh Bhalla\corref{mycorrespondingauthor}}
\ead{asbhalla@sdsu.edu}

\address[SDSU]{Department of Mechanical Engineering, San Diego State University, San Diego, CA}
%\address[Northwestern]{Department of Engineering Sciences and Applied Mathematics, Northwestern University, Evanston, IL}
\cortext[mycorrespondingauthor]{Corresponding author}

%\begin{abstract}
%\input{Abstract}
%\end{abstract}


\end{frontmatter}

%%%%%%%%%%%%%%%%%%%%%%%%%%%%%
\section{Introduction}
\input{Introduction}
\section{Mathematical formulation}\label{sec_math_eqs}
\input{Mathematical_formulation}

%%%%%%%%%%%%%%%%%%%%%%%%%%%%%%%%%
\section{Results and discussion}
\label{sec_results_and_discussion}
\input{Results_and_discussion}

%%%%%%%%%%%%%%%%%%%%%%%%%%%%%%%%
\section{Conclusions}
\label{sec_conclusions}
\input{Conclusions}

%%%%%%%%%%%%%%%%%%%%%%%%%%%%%%%%

\section*{Acknowledgements}
R.T and A.P.S.B~acknowledge support from NSF award OAC 1931368. R.T acknowledges support from San Diego State University Graduate Fellowship award. This research is part of the Frontera computing project at the Texas Advanced Computing Center (award DMS20015).  Frontera is made possible by National Science Foundation award OAC-1818253. Computational resources provided by Fermi high performance computing cluster at San Diego State University are also acknowledged.
%%%%%%%%%%%%%%%%%%%%%%%%%%%%%%%%
\section*{Supplementary material}
\label{supp_material}
\input{Supplementary_Material}

%%%%%%%%%%%%%%%%%%%%%%%%%%%%%%%%

%\appendix
%\renewcommand\thesection{\Alph{section}}
%\input{Appendix}

%%%%%%%%%%%%%%%%%%%%%%%%%%%%%%%%%
%%%%%%%%%%%%%%%%%%%%%%%%%%%%%%%%%
\section*{Bibliography}
\begin{flushleft}
 \bibliography{Brinkman_paper_biblography}
\end{flushleft}

\end{document}

%% file: Introduction.tex
Numerical simulation of multiphysics problems within complex domains has garnered much interest in the past couple of decades. In the seminal work by Angot et al.~\cite{Angot1999}, the authors describe a simple approach for simulating
the incompressible flow over obstacles by applying an additional forcing term to the governing equations.
In~\cite{Angot1999}, this volume penalization (VP) methodology (also known as the Brinkman penalization method) was
used to impose no-slip Dirichlet boundary conditions at the obstacle interface. Due to the simplicity
of its formulation and implementation, the VP technique has been successfully applied to study a variety of fluid-structure interaction problems,
including but not limited to water entry/exit~\cite{BhallaBP2019}, wave energy conversion~\cite{Dafnakis2020, khedkar2020inertial}, aquatic locomotion~\cite{Bergmann2011,Engels2017},
fluttering instabilities~\cite{Engels2013}, and flapping flight of insects~\cite{Kolomenskiy2011,Kolomenskiy2013}.
In all of these applications, the Dirichlet boundary condition formulation of the VP method was used.
In the past few years penalization methods for Neumann and more general Robin boundary conditions have been proposed,
although the analysis of such techniques is still an active area of research~\cite{Kadoch2012,Kolomenskiy2015,Brown2014,Hardy2019,Sakurai2019}.

Kadoch et al.~\cite{Kadoch2012} extended the Dirichlet boundary condition VP formulation of Angot et al.~\cite{Angot1999}
to allow for the imposition of homogeneous Neumann boundary conditions. Independently within the context of distributed
Lagrange multipliers based fictitious domain method, Doostmohammadi et al.~\cite{Doostmohammadi2014} informally described a way to enforce
homogeneous flux boundary conditions on an interface by simply setting the thermal conductivity to zero within the obstacle.
Sakurai et al.~\cite{Sakurai2019} recently developed a flux-based VP framework for imposing inhomogeneous,
spatially constant Neumann boundary conditions on the boundary of a penalization region, which formally extended the
methodology of
Kadoch et al.~\cite{Kadoch2012}. This extension enables the simulation of more complex problems within the VP framework, such as flux-driven thermal convection in irregular domains. In the flux-based VP approach of Sakurai et al., the diffusion coefficient
of the governing equation is modified and an additional forcing term is applied near the interface in order to impose the
desired flux value on the boundary. This provides a simple and efficient way to impose flux boundary conditions on embedded interfaces.

Through empirical testing of the penalized Poisson equation, Sakurai et al.~\cite{Sakurai2019} conclude that their method degrades to first-order accuracy if the embedded interface is not grid-aligned/grid-conforming despite the use of
second-order finite differences. They also conclude that if two interfaces are considered, grid-aligned or otherwise, and a \emph{different} flux boundary condition is imposed on each of them, then the method also degrades to first-order spatial accuracy. However, the method is second-order accurate for grid-aligned interfaces if
the same (spatially constant) Neumann boundary condition values are considered.

In this letter, we provide counter-examples to demonstrate that it is possible to retain second-order accuracy using Sakurai et al.'s method, even when different flux boundary conditions are imposed on multiple interfaces that do not conform to the Cartesian grid. We consider both continuous and discontinuous indicator functions in our test problems.  Both indicator functions yield a similar convergence rate for the problems considered here.  We also find that the order of accuracy results for some of the cases presented in Sakurai et al.  are not \emph{reproducible}. This is demonstrated by re-considering the same one- and two-dimensional Poisson problems solved in~\cite{Sakurai2019} in this letter.

The results shown in this letter demonstrate that the spatial order of accuracy of the flux-based VP approach of Sakurai et al. is between $\mathcal{O}$(1) and $\mathcal{O}$(2), and it depends on the underlying problem/model. The spatial order of accuracy cannot simply be deduced \emph{a priori} based on the imposed flux values, shapes, or grid-conformity of the interfaces, as concluded in Sakurai et al.~\cite{Sakurai2019}.  Further analysis is required to understand the spatial convergence rate of the flux-based VP method.

%% file: Mathematical_formulation.tex
Consider the computational domain $\Omega = \Omegaf \cup \Omegas$ consisting of disjoint fluid and solid
regions $\Omegaf$ and $\Omegas$, respectively. As described by Sakurai et al.~\cite{Sakurai2019},
the volume penalized Poisson equation with Neumann boundary conditions
$\ndot \grad q = \qn$ imposed on $\partial \Omegas$ reads,
\begin{equation}
-\div \left[\left\{\kappa \left(1 - \chi\right) + \eta \chi \right\} \grad q \right] = \left(1 - \chi\right) f + \div \left(\chi \vbeta\right) - \chi \div \vbeta, \label{eqn_vp_poisson}
\end{equation}
in which $q(\x)$ is a scalar quantity of interest, $\kappa(\x)$ is a spatially varying diffusion coefficient,
$f(\x)$ is a general forcing function, $\eta > 0$ is the penalization parameter,
$\fb(\x) =  \div \left(\chi \vbeta\right) - \chi \div \vbeta$ is an additional forcing function required to impose Neumann boundary
conditions on $\partial \Omegas$, and $\chi(\x)$ is an indicator
function that is $1$ if $\x \in \Omegas$ and $0$ if $\x \in \Omegaf$. The vector-valued flux forcing function $\vbeta(\x)$ is
selected such that $\ndot \vbeta = \qn$ on the interface $\partial \Omegas$. The unit normal vector $\n$ points out from
the fluid region and into the solid region. Note that since $\chi(\x) = 0$ within the
fluid domain, Eq.~\eqref{eqn_vp_poisson} simplifies to the standard Poisson equation in $\Omegaf$
\begin{equation}
-\div \left[\kappa \grad q\right] =  f.
\label{eqn_poisson}
\end{equation}

Sakurai et al.'s volume penalization formulation can also be generalized to handle problems with multiple interfaces. Consider
a computational domain composed of disjoint volumetric regions $\Omega^{\text{d}}_i$ (for $i = 1, \dots, D$) and $\Omega^{\text{n}}_j$ (for $j = 1, \dots, N)$
with imposed Dirichlet and Neumann boundary conditions, respectively;
i.e. $\Omegas = \Omega^{\text{d}}_1 \cup \Omega^{\text{d}}_2 \cup \cdots \cup \Omega^{\text{d}}_D \cup \Omega^{\text{n}}_1\cup \Omega^{\text{n}}_2 \cup \cdots \cup \Omega^{\text{n}}_N$.
The general form for the volume penalized Poisson equation is given by
\begin{equation}
-\div \left[\left\{\kappa \left(1 - \sum_{j = 1}^{N} \chi^{\text{n}}_j\right) + \sum_{j = 1}^{N} \eta \chi^{\text{n}}_j \right\} \grad q \right] =
\left(1 - \sum_{j = 1}^{N}  \chi^{\text{n}}_j \right)f + \sum_{j = 1}^{N} \left\{\div \left(\chi^{\text{n}}_j \vbeta_j \right) - \chi^{\text{n}}_j \div \vbeta_j\right\} -
\sum_{i = 1}^{D} \frac{\chi^{\text{d}}_i \left(q - \qdi\right)}{\eta}
, \label{eqn_general_vp_poisson}
\end{equation}
in which a Dirichlet boundary condition $q = \qdi$ is satisfied on $\partial \Omega^{\text{d}}_i$ and a Neumann boundary condition
$\n_j \cdot \grad q = \qnj$ is satisfied on $\partial \Omega^{\text{n}}_j$. The indicator function $\chi^{\text{n}}(\x) (\text{respectively}, \chi^{\text{d}} (\x))$ is $1$ if
$\x \in \Omega^{\text{n}} (\text{respectively}, \Omega^{\text{d}})$ and $0$ if $\x \in \Omega \setminus \Omega^{\text{n}} (\text{respectively}, \Omega^{\text{d}})$. Again, the vector-valued flux forcing functions are chosen such that
$\n_j \cdot \vbeta_j = \qnj$. Eq.~\eqref{eqn_general_vp_poisson} assumes that the same value of the penalization coefficient $\eta$ is used
for all the interfaces, although this is not an inherent limitation of this formulation nor our implementation. Note that general expressions
could also be written for other governing equations such as the heat, advection-diffusion, and incompressible Navier-Stokes equations,
however, we omit them for brevity.

Sakurai et al. considered a spatially constant value for $\qn$ (and $\qd$) along an interface embedded in a Cartesian domain and used periodic boundary conditions on $\partial \Omega$ in their work. We remark that the VP form of Eq.~\eqref{eqn_general_vp_poisson} is also valid for spatially varying $\qn(\x)$ or $\qd(\x)$, and the method is equally applicable when non-periodic boundary conditions are imposed on $\partial \Omega$. The discretization of Eqs.~\eqref{eqn_poisson}
and~\eqref{eqn_general_vp_poisson} uses standard second-order finite differences on a Cartesian grid. Moreover, the level set methodology
is used for representing the embedded interface, and standard regularized Heaviside functions are used to compute the smooth indicator function
\begin{align}
\chi(\x)&=
\begin{cases}
       1,  & \phi(\x) < -\ncells\; h,\\
       1 - \frac{1}{2}\left(1 + \frac{1}{\ncells h} \phi(\x)+ \frac{1}{\pi} \sin\left(\frac{\pi}{ \ncells h} \phi(\x)\right)\right) ,  & |\phi(\x)| \le \ncells\; h,\\
        0,  & \textrm{otherwise},
\end{cases}       \label{eqn_chi}
\end{align}
in which $\phi(\x)$ is a signed distance function to the interface that is negative when $\x \in \Omegas$ and positive when
$\x \in \Omegaf$, $h$ is the uniform grid spacing for the Cartesian mesh, and $\ncells \in \mathbb{R}$ is the number of cells over which the indicator
function is smoothed on either side of the interface. Sakurai et al. considered a discontinuous indicator function in their test problems, which we write below
\begin{align}
\chi(\x)&=
\begin{cases}
       1,  & \phi(\x) < 0,\\
       \half,  & \phi(\x) = 0,\\
        0,  & \textrm{otherwise}.
\end{cases}       \label{eqn_discontinuous_chi}
\end{align}
To compare our results with those reported in~\cite{Sakurai2019}, we also consider the discontinuous indicator function (along with the continuous one) in our test problems, .

The VP method described in this section is implemented within the open-source IBAMR software package~\cite{IBAMR-web-page}. 
We refer interested readers to our previous work for a more detailed discussion on the discretization techniques and linear solvers used within IBAMR~\cite{Nangia2019MF,Nangia2019WSI,BhallaBP2019}.
%Notes:
% Level set methodology is used to track the interface

%% file: Results_and_discussion.tex
%%%%%%%%%%%%%%%%%%%%%%%%%%%%%%%%%%%%%%%%%%%%%%%%%%%%%%%%%%%%%%%%%%%%%%%%%%%%%%%%%%%%%%%%%
In this section, we use the method of manufactured solution (MMS) to assess the accuracy of the flux-based VP approach introduced in Sakurai et al. using several examples. We discretely solve the VP Poisson Eq.~\eqref{eqn_vp_poisson}, which yields a numerical solution that approximates
$\qexact(\x)$ in the fluid domain $\Omegaf$ with the desired boundary conditions imposed on $\partial \Omegas$.
In all the cases considered here, we set $\kappa = 1$.
Eq.~\eqref{eqn_vp_poisson} is solved over the computational domain $\Omega$ and the domain is discretized with $N$
 and $N \times N$ Cartesian grid cells for 1D and 2D problems, respectively. A solid region $\Omegas$ is embedded within $\Omega$, and the grid does not conform
to its boundary.
Inhomogeneous Neumann boundary conditions $\n \cdot \nabla q =  \n \cdot \nabla \qexact = \n \cdot \vbeta$
are imposed on the boundary of the solid region $\partial \Omegas$. % Note that $\vbeta = \cos(x)\sin(y)\, \uvec{\i} + \sin(x)\cos(y)\, \uvec{\j}$
% is a spatially varying quantity in this work, while only spatially constant values of $\vbeta$ were considered in~\cite{Sakurai2019}.
As discussed by Sakurai et al.~\cite{Sakurai2019}, it is not necessary that $\vbeta = \grad \qexact$: an arbitrary function $\vbeta = \grad \widetilde{q}$ can also be used as long as $\grad \widetilde{q} = \grad \qexact$ on $\partial \Omegas$. Indeed, this would be the case in practice, as the solution to the Poisson equation is sought and not known $\mathit{a\; priori}$. The number of interface cells $\ncells$ is set to $2$ for the continuous indicator function unless otherwise stated. The order of accuracy results presented here are
determined based on the $L^1$ and $L^\infty$ norm of the error (denoted $\mathcal{E}_1$ and $\mathcal{E}_\infty$, respectively)
between the analytical and numerical solutions, which are computed only in the fluid domain. The penalization parameter
$\eta$ is chosen to be $\eta = 10^{-8}$, which is the
penalization value specified in~\cite{Sakurai2019}.

%In what follows, we first present the order of accuracy of the solution by considering the same one- and two-dimensional Poisson problems that were solved in~\cite{Sakurai2019}; we obtain different results compared to those reported in Sakurai et al.~\cite{Sakurai2019}. Next, we consider the same one-dimensional Poisson problem of Sakurai et al., but solve it on a slightly different fictitious domain; the order of accuracy of the solution is different compared to that obtained using the original domain.  We then assess the accuracy of the solution by considering geometrically complex interfaces and solve a set of two-dimensional Poisson problems; we obtain second-order accuracy for a particular choice of manufactured solution in all these cases. Finally, we present the order of accuracy of the solution by coupling the scalar transport equation with an incompressible Navier-Stokes fluid solver; second-order accuracy of the solution is preserved.

%%%%%%%%%%%%%%%%%%%%%%%%%%%%%%%%%%%%%%%%%%%%%%%%%%%%%%%%%%%%%

\subsection{Analysis of 1D Poisson equation with same inhomogeneous Neumann boundary condition} \label{sec_1d_same_flux_bc}
\begin{figure}[]
\centering
\includegraphics[scale = 0.08]{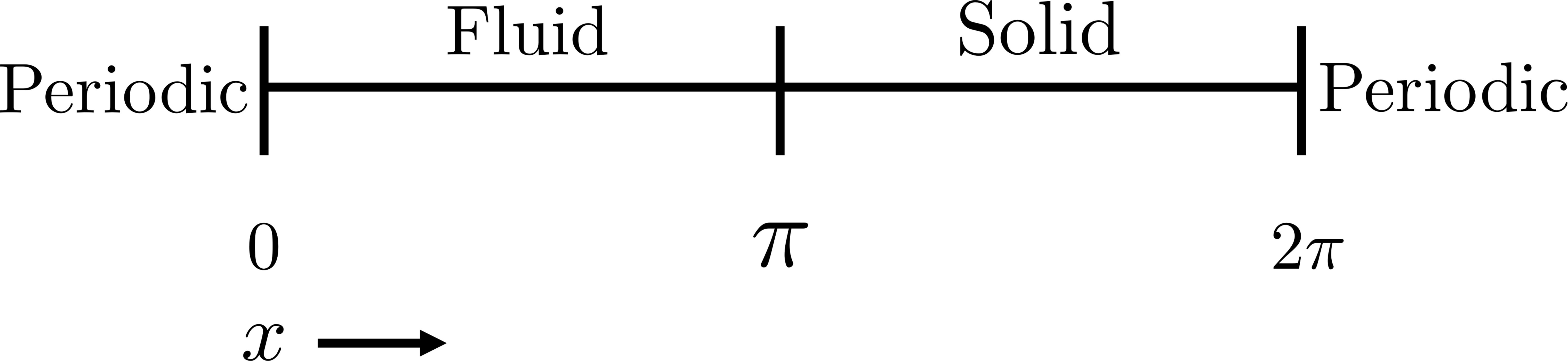}
  \caption{Schematic of the computational domain used in Sakurai et al.~\cite{Sakurai2019} to solve the 1D Poisson problem in the fluid region. The solid region in the figure represents the fictitious domain. Periodic boundary conditions are imposed on the external boundaries.}
  \label{fig_sakurai_schematic}
\end{figure}

We first consider the 1D Poisson problem with the same flux boundary condition on the two ends of the fluid domain $\Omegaf$, as done in Sec. 2.1 of Sakurai et al.~\cite{Sakurai2019}. The fluid domain $\Omegaf \in [0, \pi]$ is embedded into a larger computational domain $\Omega \in [0, 2\pi]$, as shown in Fig.~\ref{fig_sakurai_schematic}. Same inhomogeneous Neumann boundary condition value is imposed on the two fluid-solid interfaces located at $x = 0$ and $x = \pi$, respectively, and is taken to be
\begin{equation}
\left.\frac{\mathrm{d} q}{\mathrm{d} x}\right|_{x = 0} = \alpha \hspace{1 pc} \text{and} \hspace{1 pc} \left.\frac{\mathrm{d} q}{\mathrm{d} x}\right|_{x = \pi} = \alpha.
\end{equation}
We take the flux forcing function to be $\vbeta = \alpha \; \hat{\V i}$ for this test case. Here, $\hat{\V i}$ denotes the unit vector in the positive $x-$direction. The forcing function $f(\x)$ is taken to be
\begin{equation}
f(x) = m^2\cos(mx).
\end{equation}
The analytical solution of this problem using a zero-mean condition on $q$ in $\Omega_f$,  $\int _{\Omegaf} q(x) \; \text{d}x = 0$, reads as
\begin{equation}
\qexact(x) = \cos(mx) + \alpha x - \frac{\pi \alpha} {2}.
\end{equation}

 We solve the penalized Poisson equation using both continuous and discontinuous indicator functions as defined in Sec.~\ref{sec_math_eqs}.  The parameters $\alpha$ and $m$ are taken to be $1$, and periodic boundary conditions are imposed on $\partial \Omega$ (see Fig.~\ref{fig_sakurai_schematic}). Since the solution to the Poisson equation on a periodic computational domain is determinable only up to an additive constant, the discrete set of equations for this case results in a singular matrix. To invert the matrix using a direct solver, we replace the first linear equation with $\int _{\Omega} q(x) \; \text{d}x = 0$ condition, as done in Kolomenskiy et al.~\cite{Kolomenskiy2015}. We remark that although the obtained numerical solution depends on the linear equation that is replaced by the zero-mean condition (as also noted in~\cite{Kolomenskiy2015}), the order of accuracy of the solution remains the same.

\begin{figure}[]
\centering
 \subfigure[Convergence rate using interface conforming grid]{
\includegraphics[scale = 0.08]{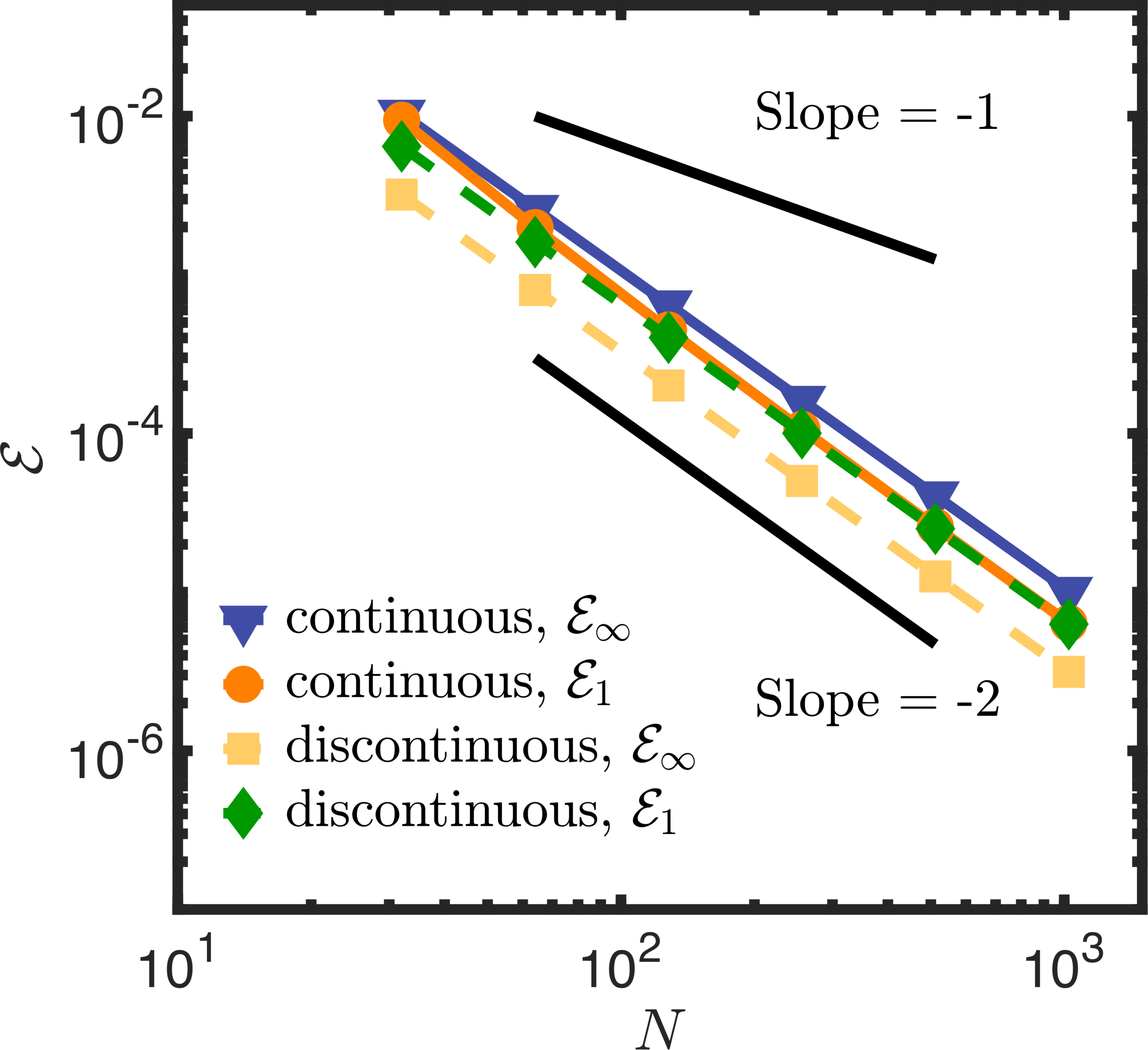}
\label{fig_sakurai_same_flux_order_of_convergence_aligned}
}
\subfigure[Convergence rate using interface non-conforming grid]{
\includegraphics[scale = 0.08]{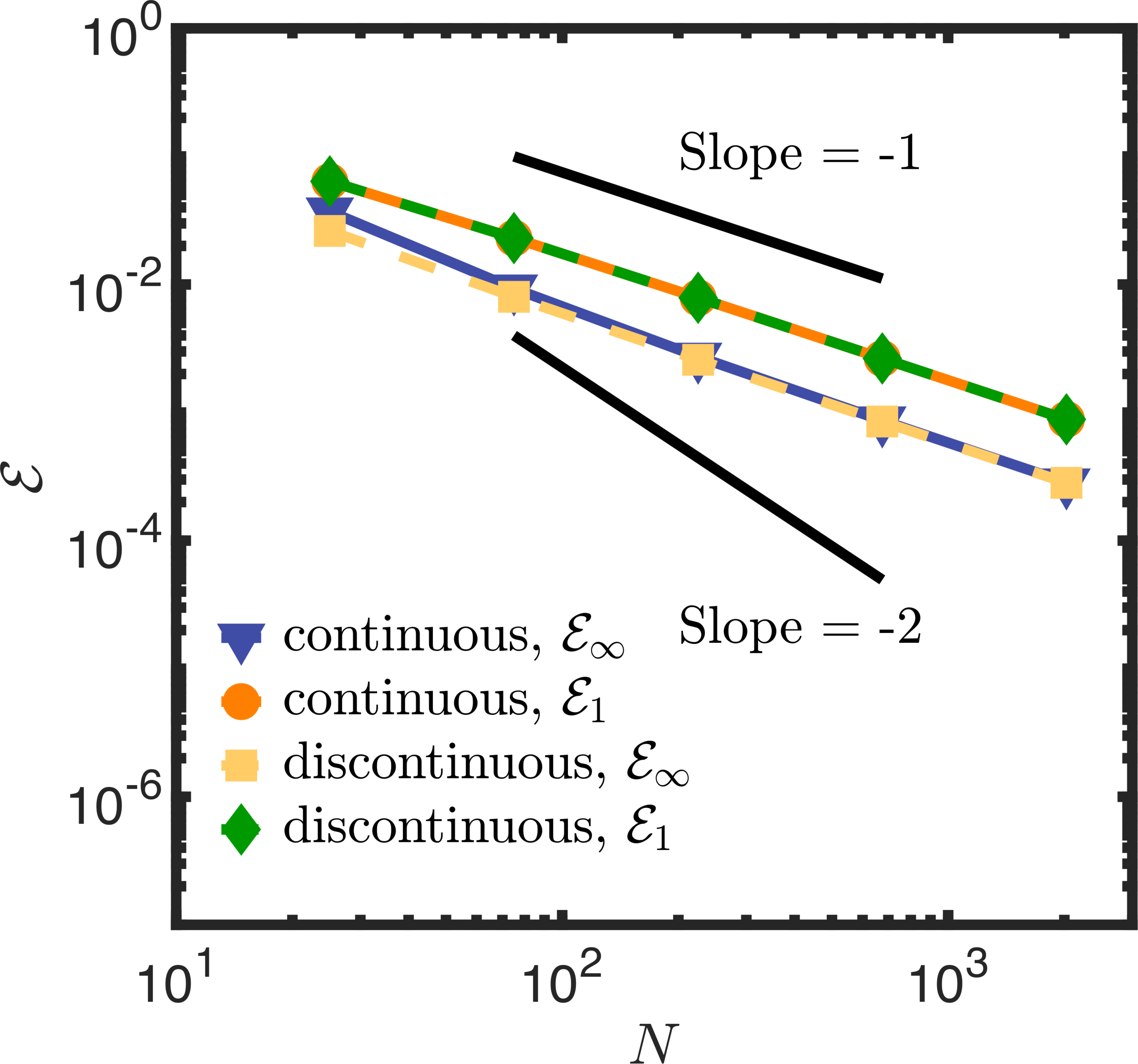}
\label{fig_sakurai_same_flux_order_of_convergence_non_aligned}
}
\subfigure[{Solution}]{
\includegraphics[scale = 0.08]{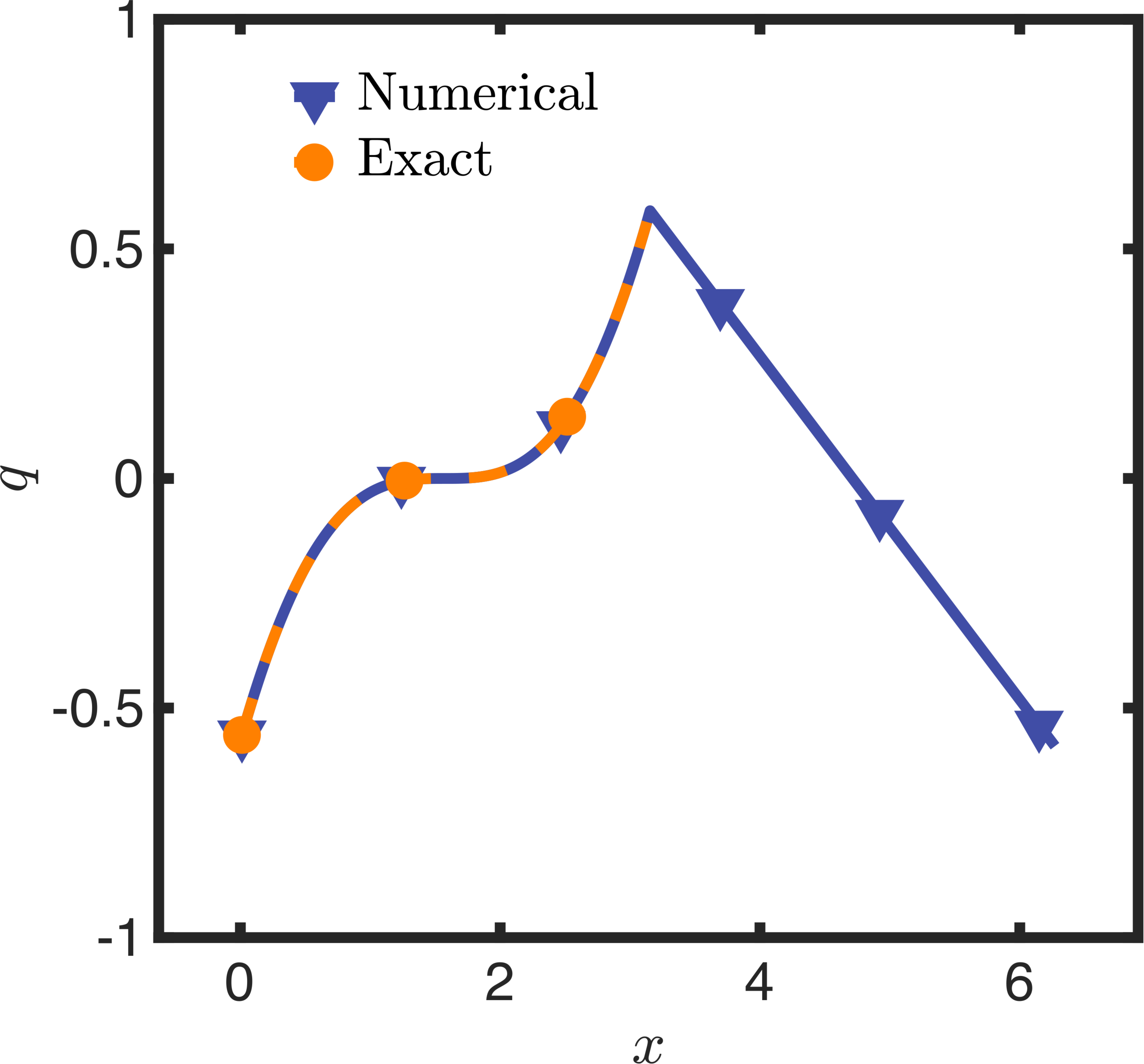}
\label{fig_sakurai_same_flux_solution}
}
  \caption{1D Poisson problem with same flux boundary conditions. Error norms $\mathcal{E}_{1}$ and $\mathcal{E}_\infty$ as a function of grid size $N$ using continuous (solid line with symbols) and discontinuous (dashed line with symbols) indicator functions when~\subref{fig_sakurai_same_flux_order_of_convergence_aligned} the fluid-solid interface at $x = \pi$ is aligned with the Cartesian cell face, and~\subref{fig_sakurai_same_flux_order_of_convergence_non_aligned} when it is not; \subref{fig_sakurai_same_flux_solution} numerical solution $q$ obtained using $N = 256, \alpha = 1$, and $m =1$, along with the exact solution. The penalization parameter $\eta$ is taken as $10^{-8}$.}
\label{fig_sakurai_same_flux}
\end{figure}

Two sets of $N$ values are selected to assess the order of accuracy of the solution: (i) $\Nalign$ = [32, 64, 128, 256, 512,1024] which aligns the fluid-solid interface
located at $x = \pi$ with the Cartesian cell face, as  done in~\cite{Sakurai2019}, and (ii) $\Nnalign$ = [25, 75,  225, 675, 2025] which does not. The other two fluid-solid interfaces at  $x = 0$ and $x = 2\pi$ are located on grid cell faces by construction. Fig.~\ref{fig_sakurai_same_flux} compares the spatial convergence rate for the two grid setups. As observed in Fig.~\ref{fig_sakurai_same_flux_order_of_convergence_aligned}, when the interface aligns with the Cartesian grid face, $\mathcal{O}(2)$ convergence rate is obtained using both continuous  and discontinuous indicator functions. Second-order spatial accuracy is also obtained in Sakurai et al.~\cite{Sakurai2019} using the discontinuous indicator function using a similar grid setup. However, the order of accuracy degrades to  $\mathcal{O}(1)$ when the interface is not aligned with the grid, as observed in Fig.~\ref{fig_sakurai_same_flux_order_of_convergence_non_aligned}. The authors in~\cite{Sakurai2019} did not present the order of accuracy results using a non-conforming grid (to the interface) for this problem.  Finally, Fig.~\ref{fig_sakurai_same_flux_solution} shows the numerical solution $q$, and compares it against the exact solution for $N = 256$ grid. An excellent agreement is obtained.

%%%%%%%%%%%%%%%%%%%%%%%%%%%%%%%%%%%%%%%%%%%%%%%%%%%%%%%%%%%%%%%%%%%%

\subsection{Analysis of 1D Poisson equation with different inhomogeneous Neumann boundary conditions} \label{sec_1d_different_flux_bc}
\begin{figure}[]
\centering
 \subfigure[Convergence rate using interface conforming grid]{
\includegraphics[scale = 0.08]{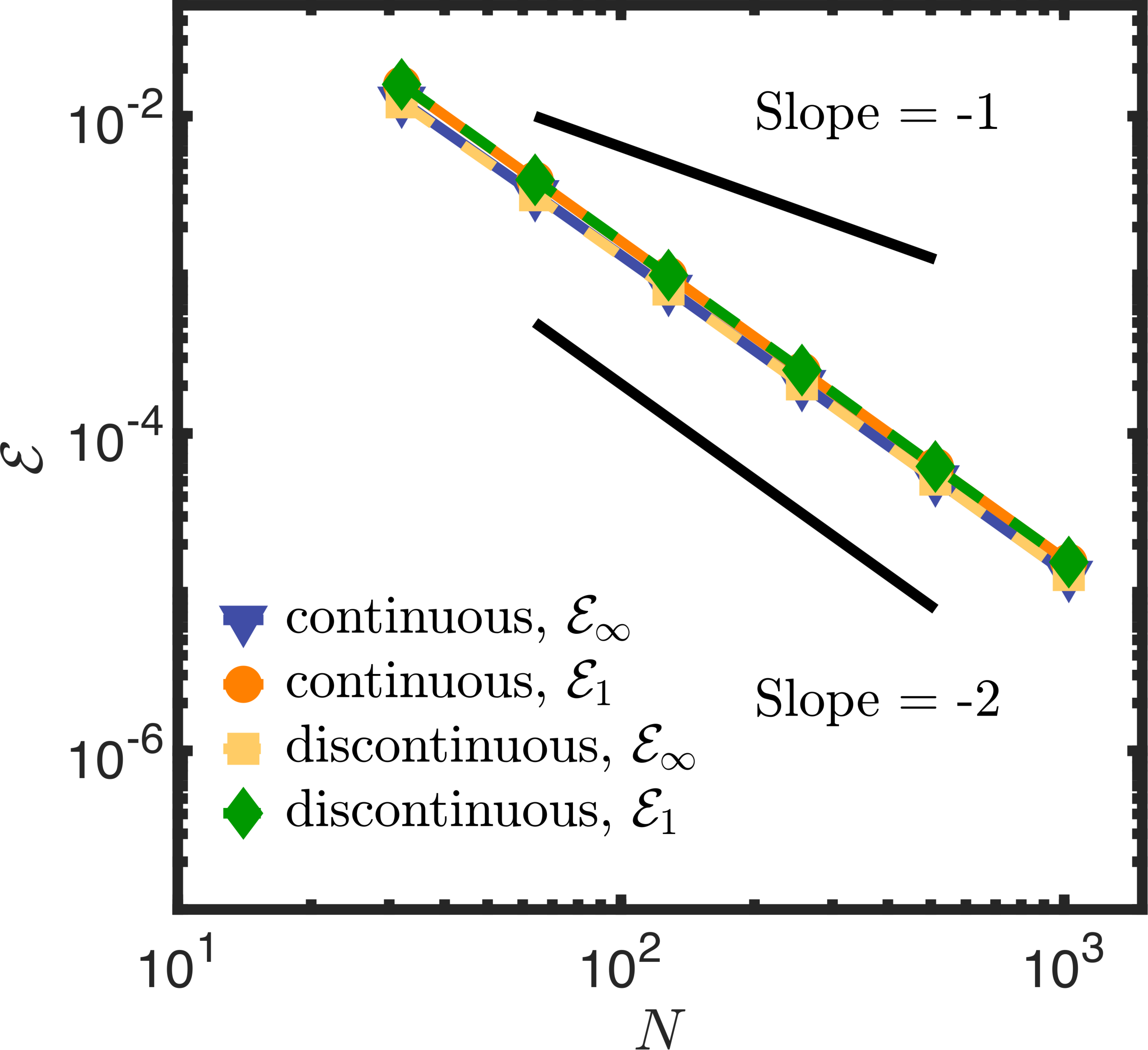}
\label{fig_sakurai_different_flux_order_of_convergence_aligned}
}
 \subfigure[Convergence rate using interface non-conforming grid]{
\includegraphics[scale = 0.08]{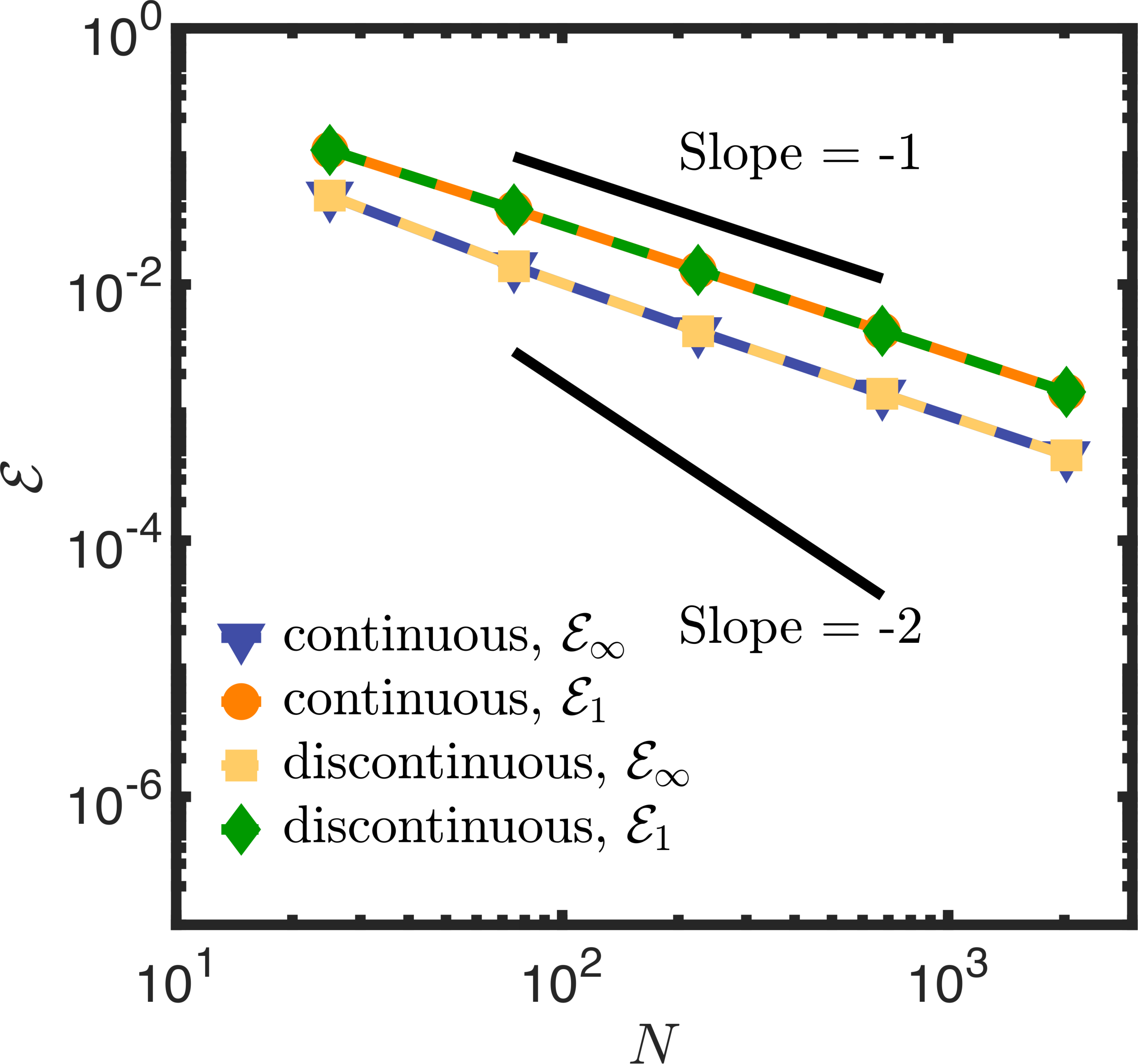}
\label{fig_sakurai_different_flux_order_of_convergence_non_aligned}
}
\subfigure[{Solution}]{
\includegraphics[scale = 0.08]{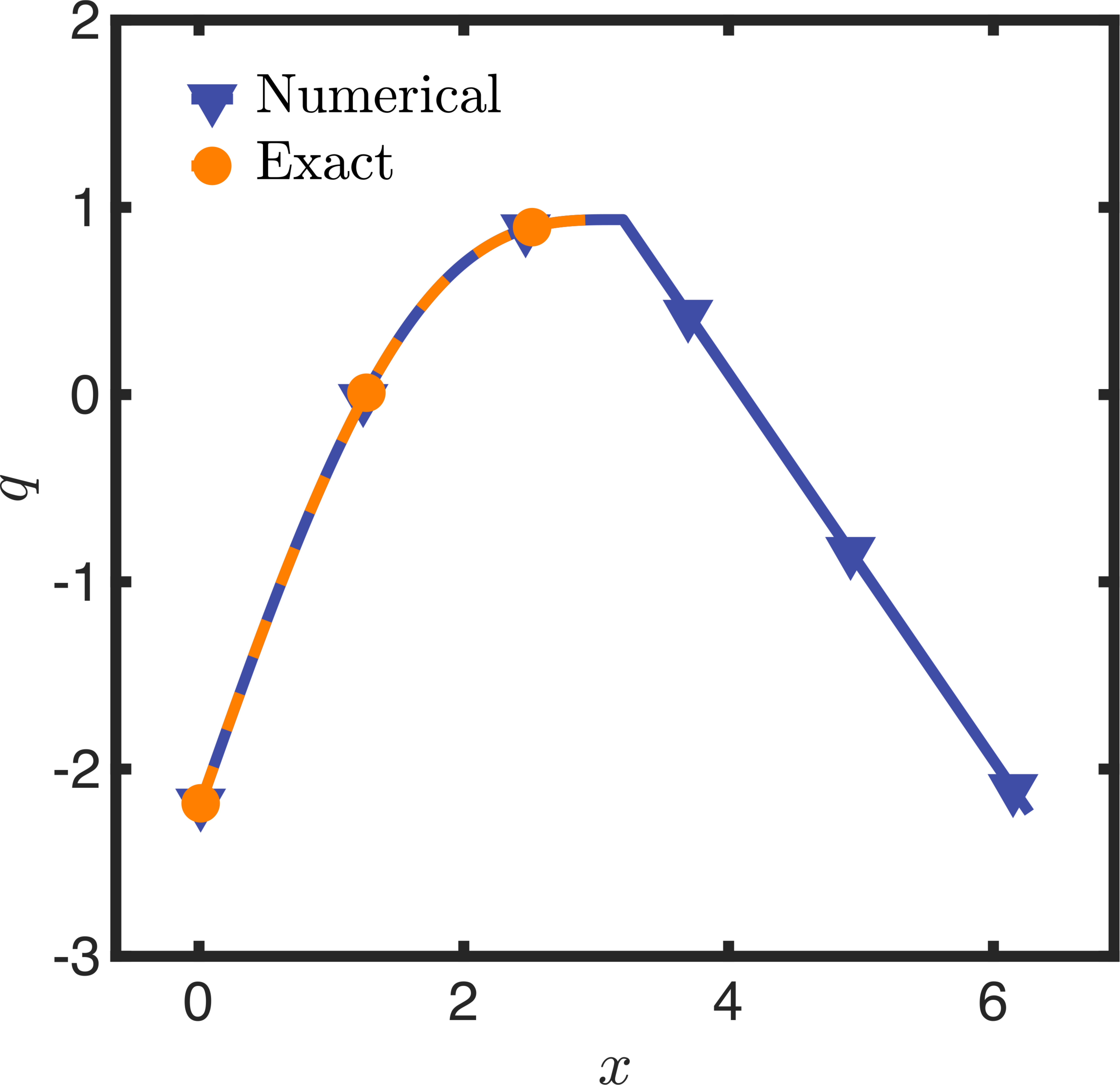}
\label{fig_sakurai_different_flux_solution}
}
  \caption{1D Poisson problem with different flux boundary conditions. Error norms $\mathcal{E}_{1}$ and $\mathcal{E}_\infty$ as a function of grid size $N$ using continuous (solid line with symbols) and discontinuous (dashed line with symbols) indicator functions when~\subref{fig_sakurai_different_flux_order_of_convergence_aligned}  the fluid-solid interface at $x = \pi$ is aligned with the Cartesian cell face, and~\subref{fig_sakurai_different_flux_order_of_convergence_non_aligned} when it is not; \subref{fig_sakurai_different_flux_solution} numerical solution $q$ obtained using $N = 256, \alpha = 1$, and $m =1$, along with the exact solution. The penalization parameter $\eta$ is taken as $10^{-8}$.}
\label{fig_sakurai_different_flux}
\end{figure}

We now consider the 1D Poisson problem with different inhomogeneous Neumann boundary conditions on the two ends of the fluid domain, as done in Sec. 2.3 of Sakurai et al. The forcing function $f(\x)$ for this case is
\begin{equation}
f(x) = m^2\sin(mx),
\end{equation}
and the inhomogeneous Neumann boundary condition values on the two ends are
\begin{equation}
\left.\frac{\mathrm{d} q}{\mathrm{d} x}\right|_{x = 0} = \alpha + m \hspace{1 pc} \text{and} \hspace{1 pc} \left.\frac{\mathrm{d} q}{\mathrm{d} x}\right|_{x = \pi} = \alpha - m.
\label{eq_1d_poisson_different_flux_bc}
\end{equation}
Here, $\alpha$ and $m$ parameters are taken to be 1. The problem setup remains the same as shown in Fig.~\ref{fig_sakurai_schematic}. The analytical solution of this problem (using a zero-mean condition on $q$ in $\Omega_f$) reads as
\begin{equation}
\qexact(x) = \sin(mx) + \alpha x - \frac{2}{m\pi} - \frac{\pi \alpha} {2}.
\end{equation}
The flux forcing function is taken to be $\vbeta = \nabla \qexact = \left(m\cos(mx) + \alpha \right)\hat{\V i}$, which also satisfies the boundary conditions written in Eq.~\eqref{eq_1d_poisson_different_flux_bc}. We again replace the first linear equation with zero mean of $q$ in $\Omega$ to obtain the unique solution. The results for this case are presented in Fig.~\ref{fig_sakurai_different_flux} for $\Nalign$ and $\Nnalign$ grid size values, as taken in the previous section.

In Fig.~\ref{fig_sakurai_different_flux_order_of_convergence_aligned} we observe $\mathcal{O}(2)$ convergence for both types of indicator functions when the interface aligns with the Cartesian grid face. This is in contrast to Sakurai et al.~\cite{Sakurai2019} where $\mathcal{O}(1)$ convergence is reported for this test problem using a similar grid setup; the results reported in~\cite{Sakurai2019} are not reproducible despite the use of same discretization method and problem setup. The authors in~\cite{Sakurai2019} attribute the reduction in accuracy to different values of flux boundary condition, which is clearly not the case here.  Fig.~\ref{fig_sakurai_different_flux_order_of_convergence_non_aligned} shows the order of accuracy results when the interface is not aligned with the Cartesian cell face --- $\mathcal{O}(1)$ convergence rate is exhibited using both continuous and discontinuous indicator functions. Finally, Fig.~\ref{fig_sakurai_different_flux_solution} shows the numerical solution $q$, and compares it against the exact solution for $N = 256$ grid. An excellent agreement is obtained.

The results presented in the above two sections may suggest that the spatial accuracy of the flux-based VP method is  $\mathcal{O}(2)$ when the interface aligns the Cartesian mesh, but degrades to  $\mathcal{O}(1)$ when it does not. This is also one of the conclusions in Sakurai et al. However, our next examples contradict this conclusion.

%%%%%%%%%%%%%%%%%%%%%%%%%%%%%%%%%%%%%%
\subsection{Analysis of 2D Poisson equation with different flux boundary conditions } \label{sec_2d_different_flux_bc}
\begin{figure}[]
\centering
\subfigure[{Schematic of the case}]{
\includegraphics[scale = 0.07]{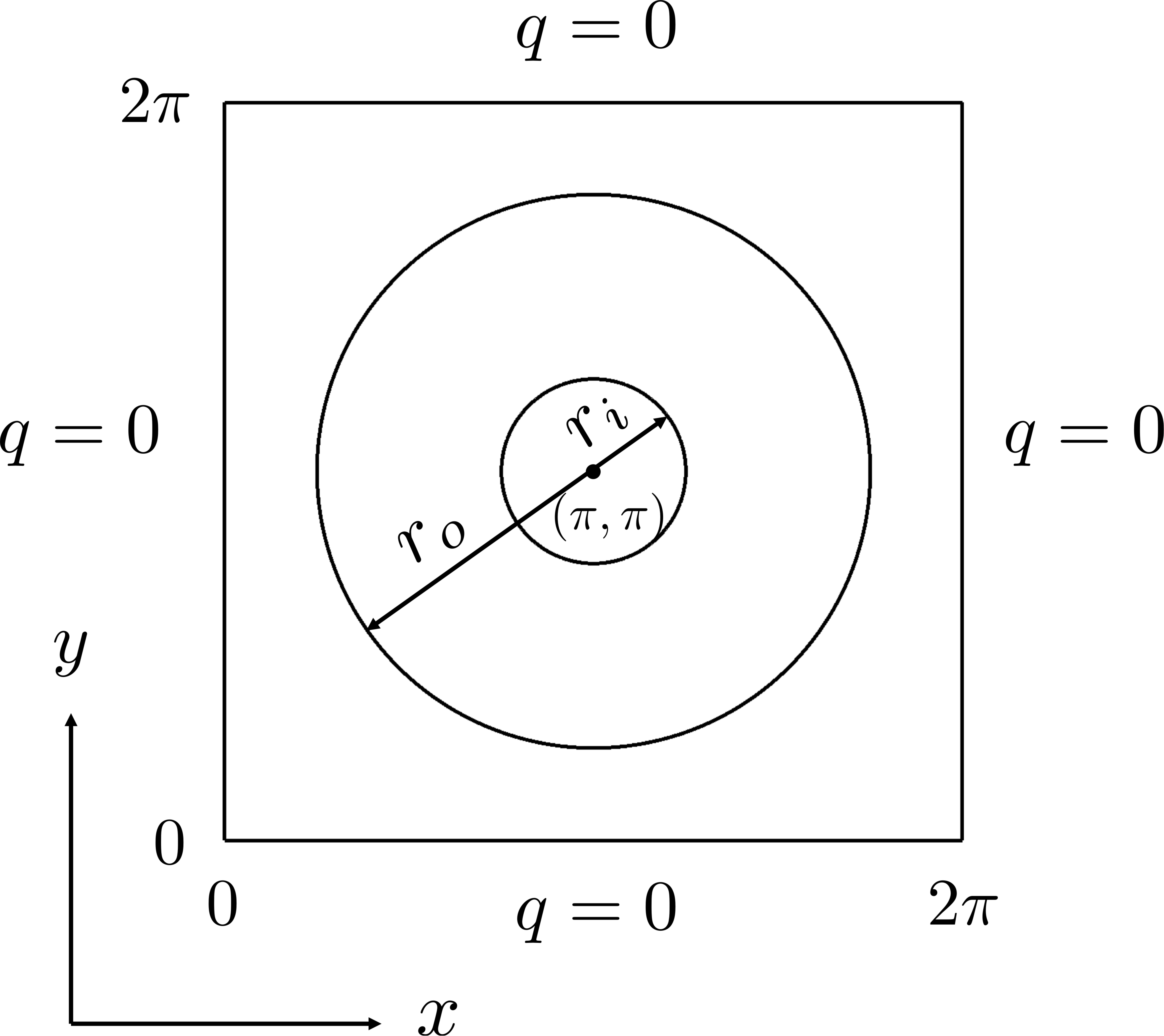}
\label{fig_sakurai_2D_Poisson_schematic}
}
\subfigure[{Solution variation along $y-$direction}]{
\includegraphics[scale = 0.08]{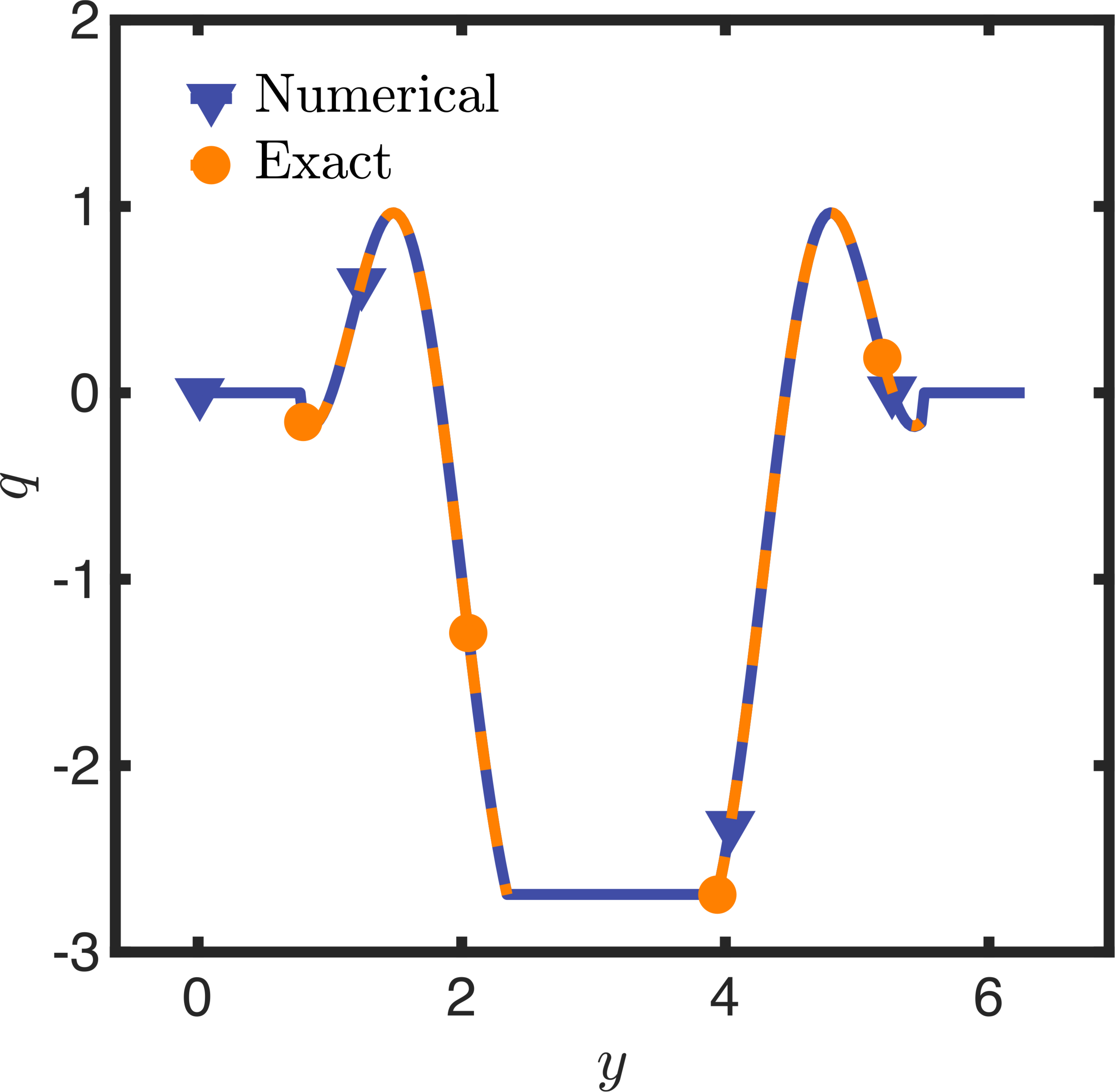}
\label{fig_sakurai_2D_Poisson_xcut_solution}
}
\subfigure[{Solution variation along $x-$direction}]{
\includegraphics[scale = 0.08]{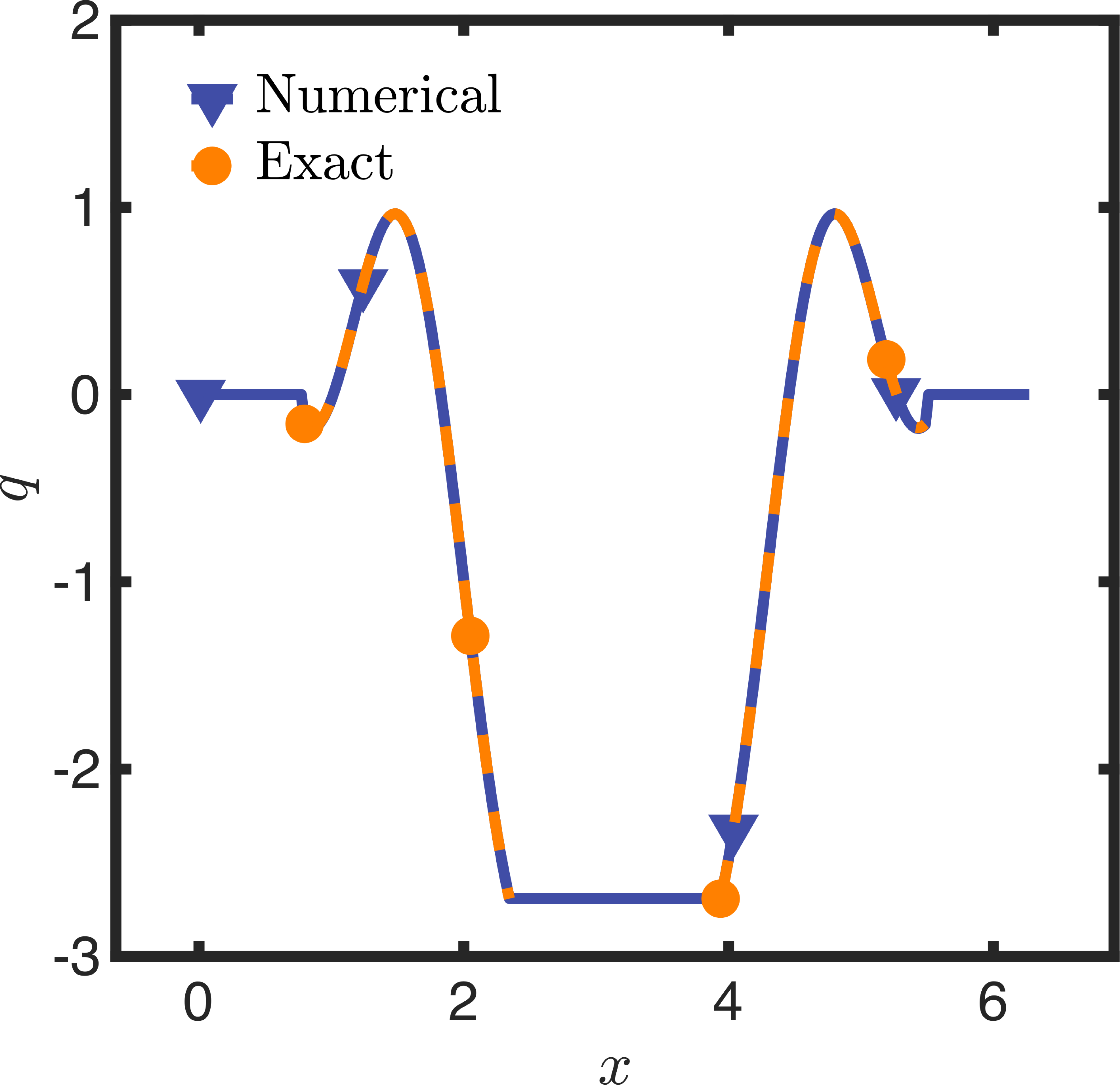}
\label{fig_sakurai_2D_Poisson_ycut_solution}
}
 \subfigure[Spatial convergence rate]{
\includegraphics[scale = 0.08]{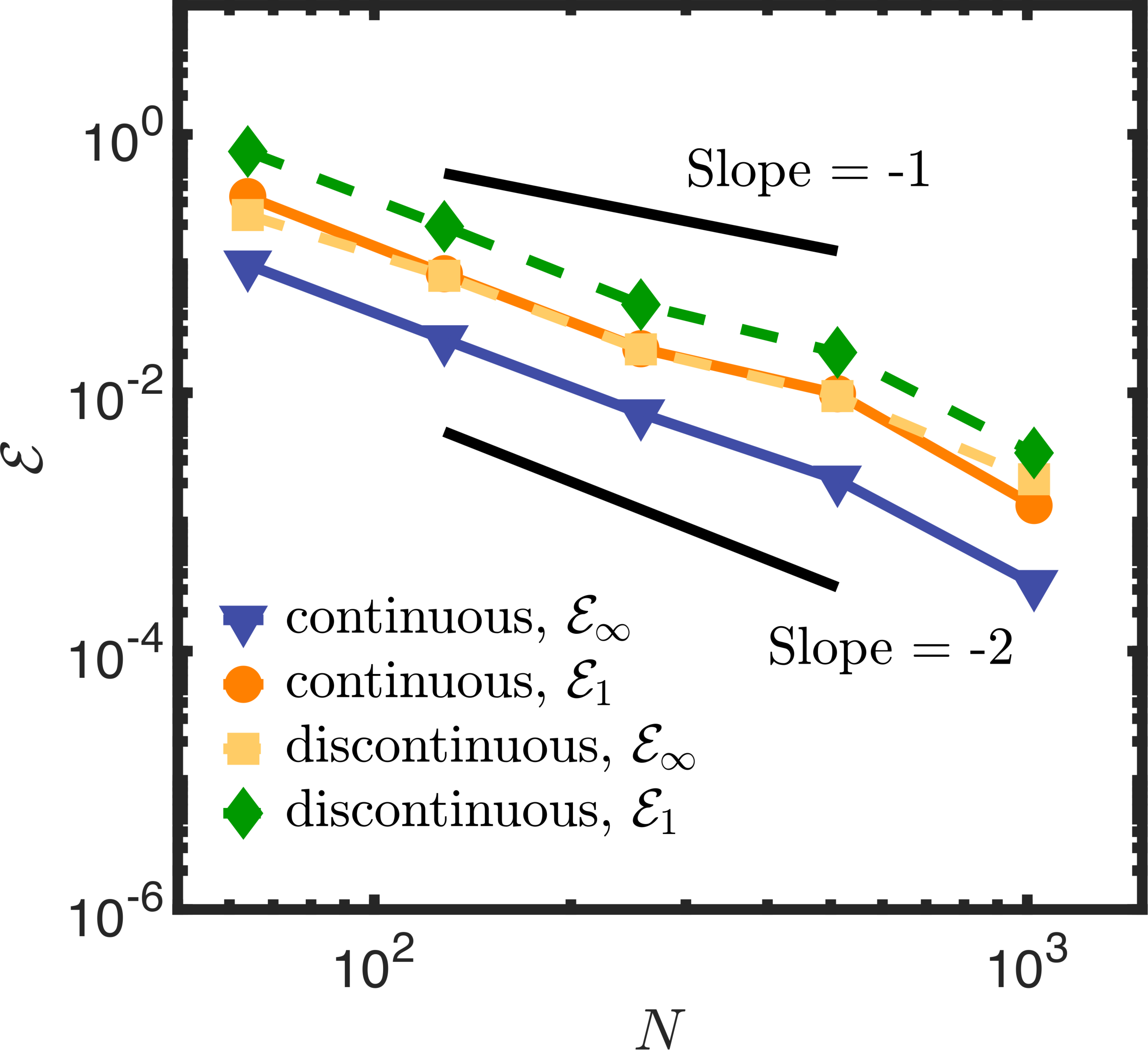}
\label{fig_sakurai_2D_Poisson_order_of_convergence}
}
 \caption{2D Poisson problem with different flux boundary conditions on the two interfaces: \subref{fig_sakurai_2D_Poisson_schematic} problem setup; \subref{fig_sakurai_2D_Poisson_xcut_solution} variation of the numerical solution along $y-$direction at a fixed $x = 3.12$ location using $N = 256$ grid;  \subref{fig_sakurai_2D_Poisson_ycut_solution} variation of numerical solution along $x-$direction at a fixed $y= 3.12$  location using $N = 256$ grid; \subref{fig_sakurai_2D_Poisson_order_of_convergence}  error norms $\mathcal{E}_{1} $ and $\mathcal{E}_\infty$ as a function of grid size $N$ using continuous (solid lines with symbols) and discontinuous (dashed lines with symbols) indicator functions. The penalization parameter $\eta$ is taken as $10^{-8}$, and  $\alpha$ is taken as 1.}
\label{fig_sakurai_2D_Poisson}
\end{figure}

Here we solve the penalized form of Poisson equation~\ref{eqn_vp_poisson} in a circular annulus using different flux boundary conditions
on the two interfaces. The same case is considered in Sec. 3 of Sakurai et al. The circular annulus has an inner radius $r_i$ of $\pi/4$ and an outer radius $r_o$ of $3\pi/4$, and is centered around the point $(\pi,\pi)$. The annulus is embedded into a larger computational domain $\Omega \in [0, 2\pi]^2$, as shown in Fig.~\ref{fig_sakurai_2D_Poisson_schematic}. The forcing function for this case is
\begin{equation}
f(r) = 16\cos(4r)+\frac{4\sin(4r)}{r},
\end{equation}
in which $r = \sqrt{\left(x - \pi\right)^2 + \left(y - \pi\right)^2}$, and the flux boundary condition values on the two interfaces are taken to be
\begin{equation}
\left.\frac{\mathrm{d} q}{\mathrm{d} r}\right|_{r = \frac{\pi}{4}} = 3\alpha \hspace{1 pc} \text{and} \hspace{1 pc} \left.\frac{\mathrm{d} q}{\mathrm{d} r}\right|_{r = \frac{3\pi}{4}} = \alpha.
\end{equation}
The exact solution for this problem using the zero-mean condition $\int_{\Omegaf}rq(r)\; \text{d}r = 0$  reads as
\begin{equation}
\qexact(r) = \cos(4r) + \frac{3}{4} \alpha \pi \log(r) - \frac{3}{32} \alpha \pi \left(9 \log\left(\frac{3}{4}\pi\right) - \log\left(\frac{\pi}{4}\right) - 4\right).
\label{sakurai_2d_poisson_exact}
\end{equation}
We solve the 2D penalized Poisson equation using both continuous and discontinuous indicator functions. Homogeneous Dirichlet boundary conditions are imposed on $\partial \Omega$. The flux forcing function $\vbeta$ for this case is taken to be  $\vbeta(\x) = g(\x)\e_r$,  in which $\e_r = \left(\frac{x-\pi}{r}, \frac{y-\pi}{r}\right)$ and $g(\x)$ is
\begin{equation}
g(\x) = \left \{  \begin{array}{ll}
      \alpha \left(\frac{4r}{3\pi}\right)^2 \left(4\left(1 - \frac{r}{\pi}\right)\right)^3 & \text{if} \hspace{1pc} 0 \leq r \leq \pi \\
      0 & \text{otherwise} \\
\end{array} \right.
\end{equation}
The same form of $\vbeta$ is also used in~\cite{Sakurai2019},  although $\vbeta = \grad \qexact$ can also be defined here. Fig.~\ref{fig_sakurai_2D_Poisson} compares the numerical solution with the exact solution as written in Eq.~\ref{sakurai_2d_poisson_exact}. The convergence rate of the solution as a function of grid size is also shown. From Fig.~\ref{fig_sakurai_2D_Poisson_order_of_convergence}, we note that the convergence rate is close to $\mathcal{O}(2)$, as opposed to  $\mathcal{O}(1)$ reported in Sakurai et al. for this problem. We again remark that despite using the same problem setup and discretization technique, the spatial order of accuracy shown in~\cite{Sakurai2019} is not reproducible for this test problem as well. Moreover, different values of the flux boundary condition and/or the circular shape of the interface did not reduce the order of accuracy, as was reasoned by the authors in~\cite{Sakurai2019} for this problem.

%%%%%%%%%%%%%%%%%%%%%%%%%%%%%%%%%%%%%%%%%%%%%%%%%%%%%%%%%%%%%%%%%%%%%%%%%%%%%%%%%%%%%%%%%

\subsection{Analysis of flux boundary condition on complex interfaces} \label{sec_complex_shapes}

In this section, we consider geometrically complex interfaces and  use a manufactured solution of the form
\begin{equation}
\qexact(\x)= \sin(x)  \sin(y),
\label{q_analytical}
\end{equation}
to demonstrate that the spatial order of accuracy for the flux-based VP method can indeed be $\mathcal{O}(2)$, despite imposing different flux values on multiple interfaces that do not conform to the Cartesian grid.

\begin{figure}[]
  \centering
   \includegraphics[scale = 0.13]{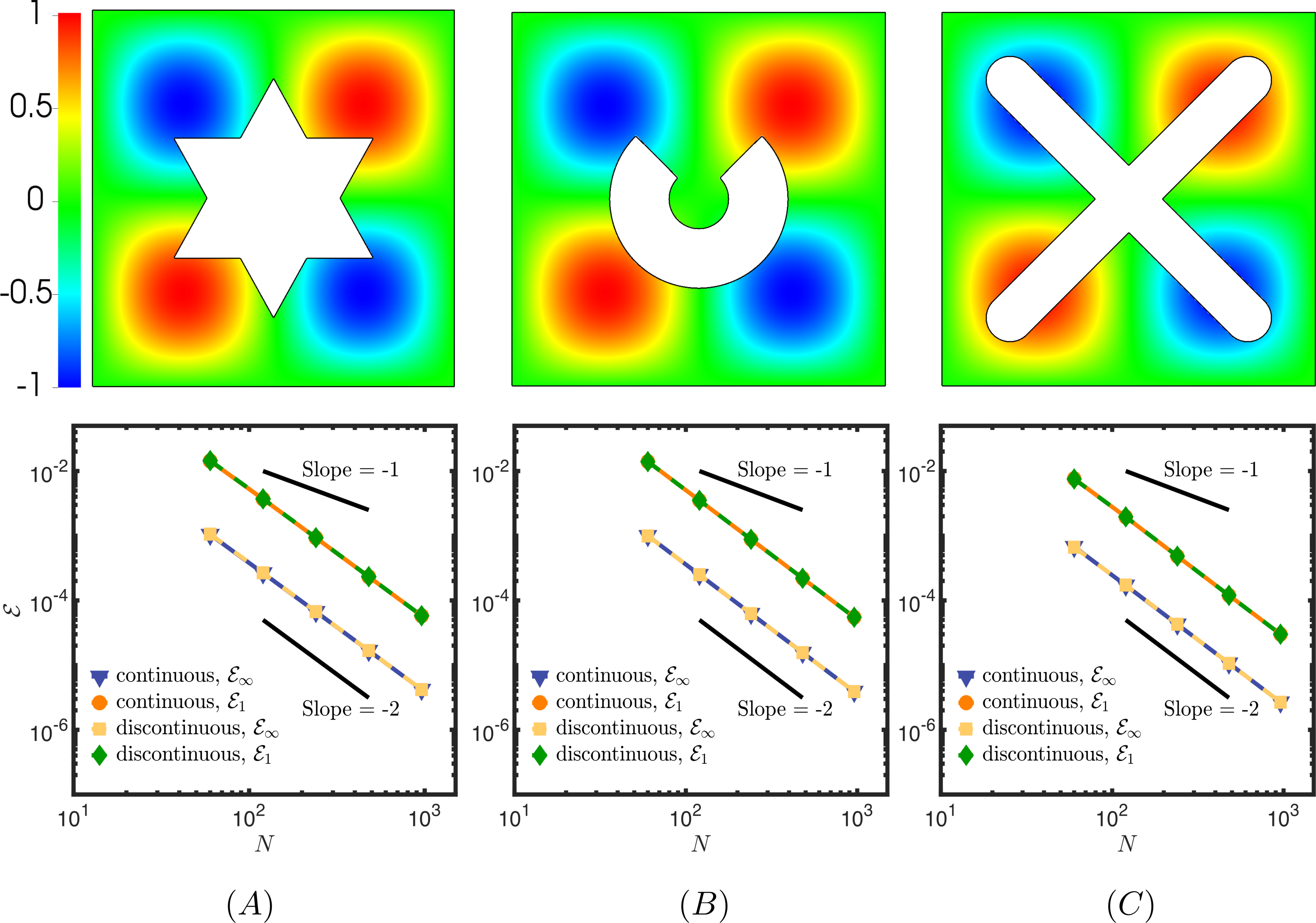}
   \caption{The numerical solution $q$ (top row), and the error norms $\mathcal{E}_{1}$ and $\mathcal{E}_\infty$ as a function of grid size $N$ (bottom row) for three complex annuli; (A) Hexagram; (B) Horseshoe; (C) X-cross. The results are shown for both continuous (solid line with symbols) and discontinuous (dashed line with symbols) indicator functions. The penalization parameter is taken as $10^{-8}$.}
     \label{fig_complex_annulus}
\end{figure}

Specifically, we consider three different interfacial geometries centered about the point$(\pi, \pi)$: a hexagram, a horseshoe, and an x-cross. The region interior to the interface is considered to be the (fictitious) solid domain; see Fig.~\ref{fig_complex_annulus}. It can be noted that these shapes involve sharp corners and the interfaces do not align with the grid.  The required forcing function $f(\x)$ is obtained by plugging Eq.~\eqref{q_analytical} into Eq.~\ref{eqn_poisson}, and the flux forcing function is taken to be $\vbeta = \grad \qexact$. Dirichlet boundary conditions are imposed on the external boundaries, i.e., $\left.q\right|_{\partial \Omega(\x)} = \qexact(\x)$, and spatially varying flux boundary conditions are imposed on the embedded interfaces.  As shown in Fig.~\ref{fig_complex_annulus}, second-order spatial convergence rates are exhibited for each of these complex annuli using both continuous and discontinuous indicator functions.
%Moreover, we note that the sharp corners seen in these cases do not further degrade the error.

\begin{figure}[]
\centering
\subfigure[{Numerical solution}]{
\includegraphics[scale = 0.14]{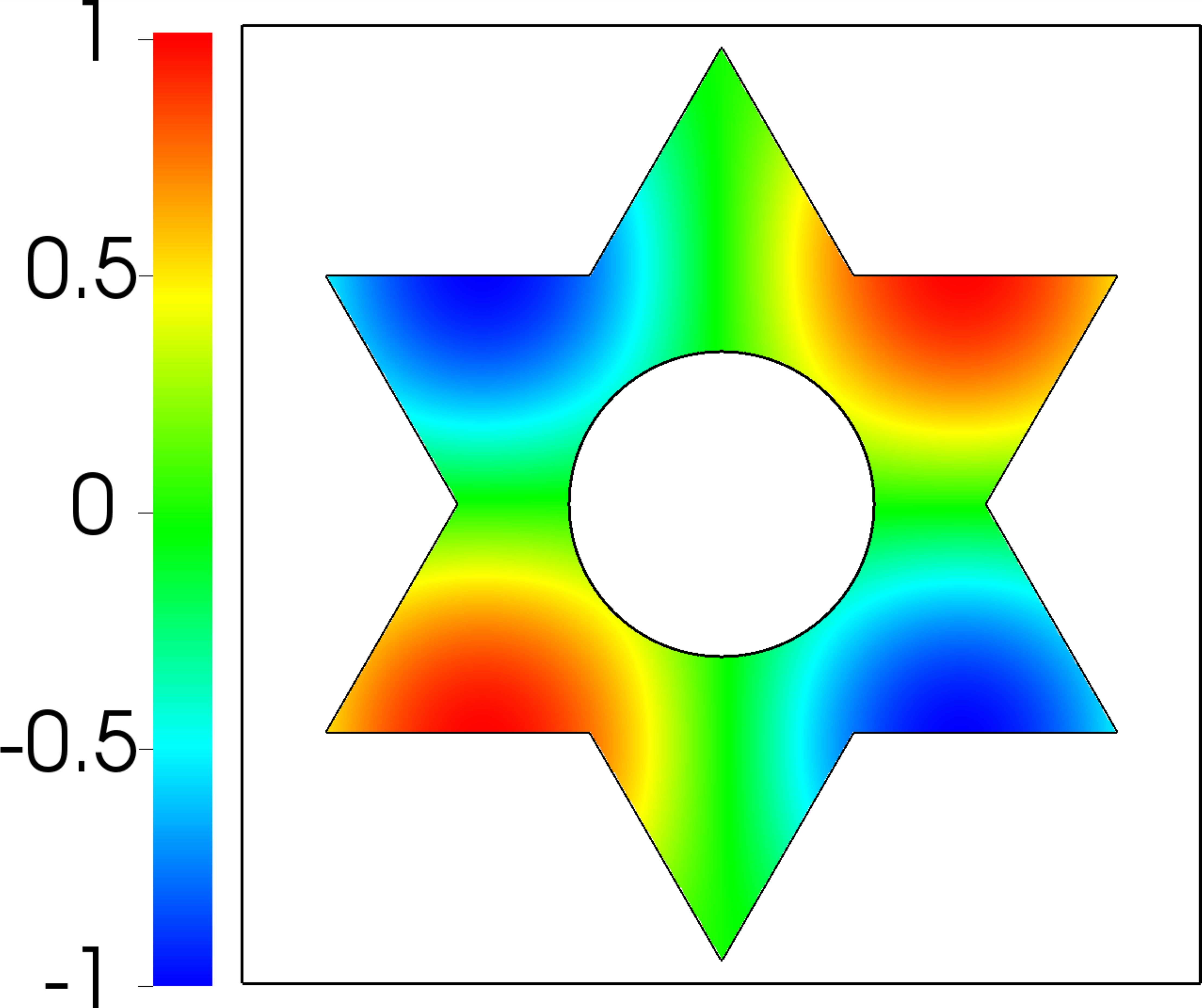}
\label{fig_hexagram_circle_solution}
}
 \subfigure[Numerical solution]{
\includegraphics[scale = 0.14]{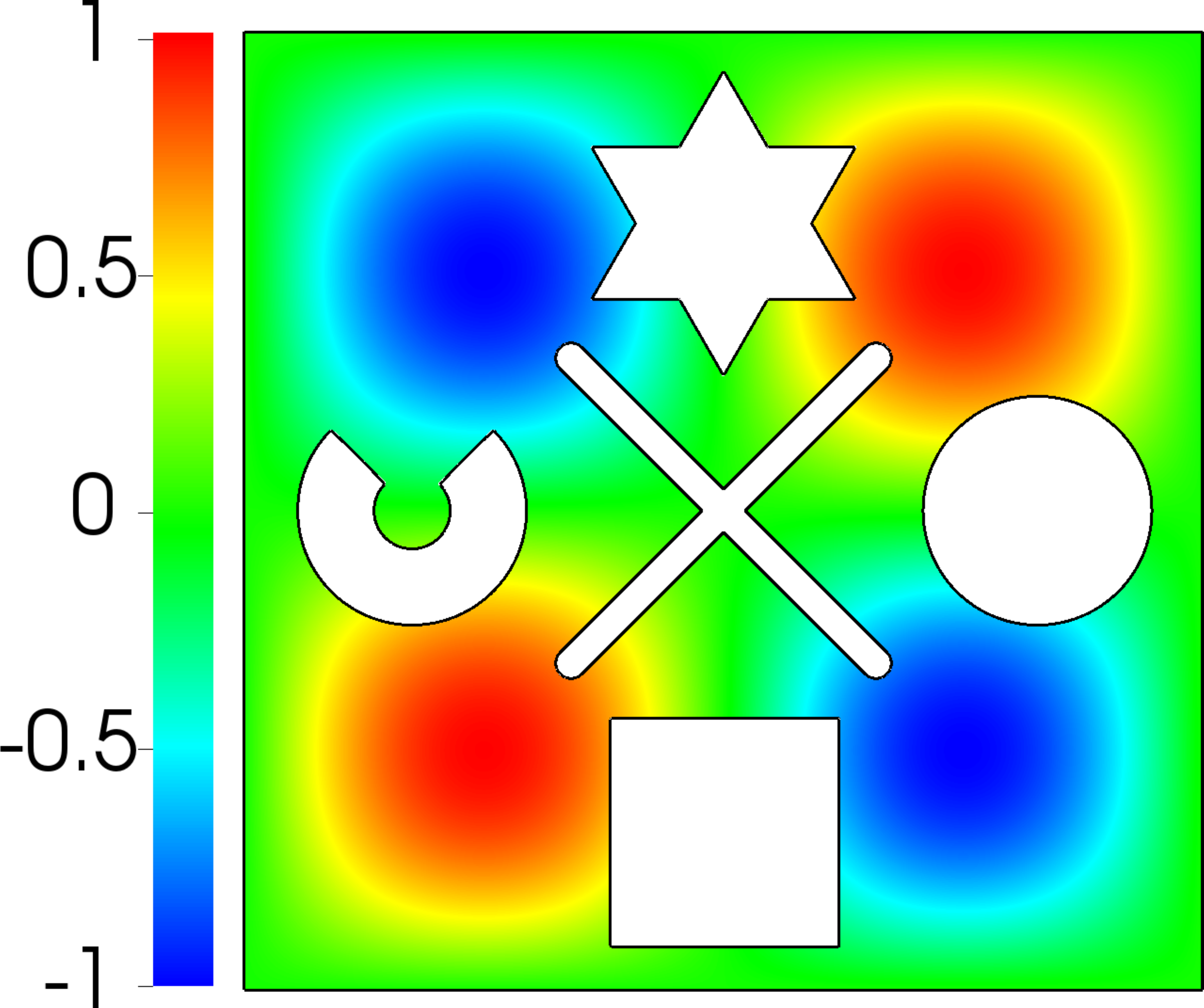}
\label{fig_all_shapes_solution}
}
\subfigure[{Spatial convergence rate}]{
\includegraphics[scale = 0.08]{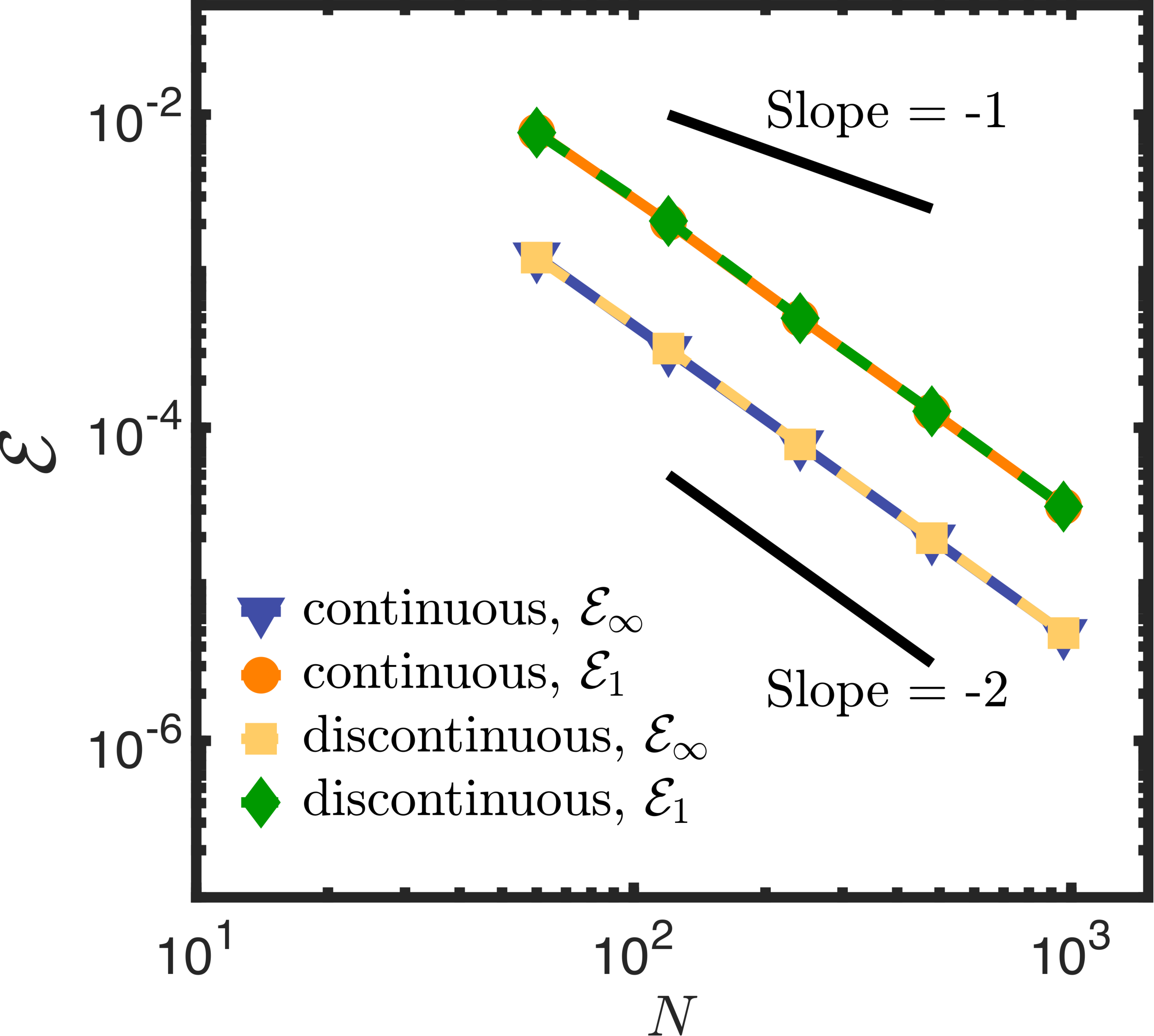}
\label{fig_hexagram_circle_order_of_convergence}
}
 \subfigure[Spatial convergence rate]{
\includegraphics[scale = 0.08]{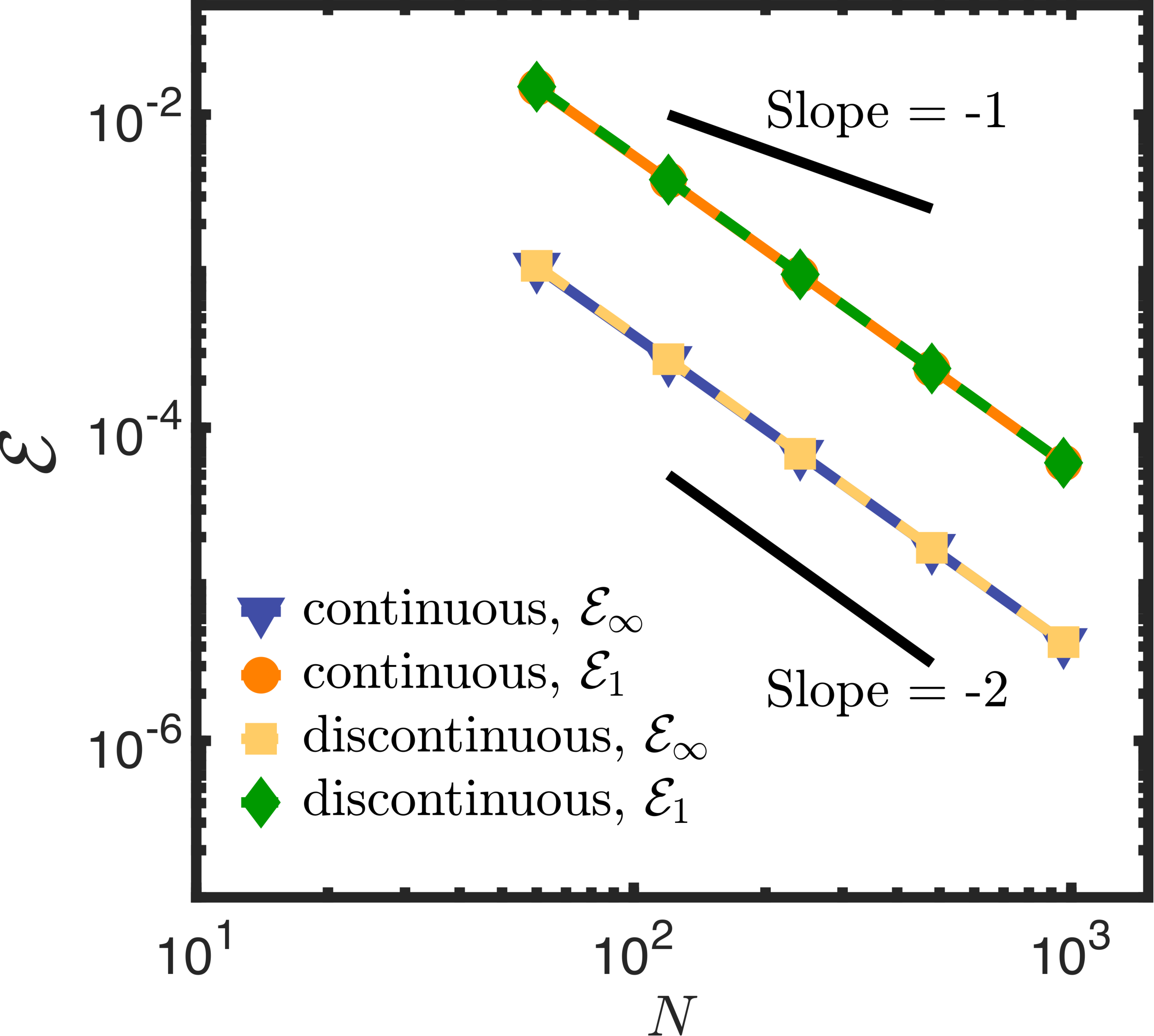}
\label{fig_all_shapes_order_of_convergence}
}

 \caption{\subref{fig_hexagram_circle_solution} Numerical solution in a complex annulus bounded by a hexagram and a circle; \subref{fig_all_shapes_solution} Numerical solution in a complex annulus formed by multiple interfaces and the rectangular computational domain;
  \subref{fig_hexagram_circle_order_of_convergence} Error norms $\mathcal{E}_{1} $ and $\mathcal{E}_\infty$ as a function of grid size $N$ for the first annulus case; \subref{fig_all_shapes_order_of_convergence} Error norms $\mathcal{E}_{1}$ and $\mathcal{E}_\infty$ as a function of grid size $N$ for the second annulus case. The results are shown for both continuous (solid line with symbols) and discontinuous (dashed line with symbols) indicator functions. The penalization parameter $\eta$ is taken as $10^{-8}$.}
\label{fig_2d_multiple_complex_interfaces}
\end{figure}

We also consider two additional complex domains with the same manufactured solution as written in Eq.~\eqref{q_analytical}. The first one is a complex annulus whose outer surface is a hexagram and the inner surface is a circle of radius 1. Both surfaces are centered about the point $(\pi, \pi)$.  Flux boundary conditions are imposed on the two surfaces of the annulus, whereas homogeneous Dirichlet boundary conditions are imposed on the external boundaries of the computational domain.  The zero-mean condition on $q$ in the fluid/annulus domain is imposed as a post-processing step to obtain the unique solution for this case. Fig.~\ref{fig_hexagram_circle_order_of_convergence} shows the spatial order of accuracy for this case. Second-order convergence is exhibited. For the second complex domain case, we embed all of the previously considered interfaces into a rectangular computational domain and impose spatially varying flux boundary condition on the interfaces. The fluid domain is exterior to all the interfaces. The penalized Poisson equation is solved by imposing $\left.q\right|_{\partial \Omega(\x)} = \qexact(\x)$. Fig.~\ref{fig_all_shapes_order_of_convergence} shows the convergence rate for the second annulus case. Again, the method exhibits $\mathcal{O}$(2) convergence.

%%%%%%%%%%%%%%%%%%%%%%%%%%%%%%%%%%%%
\subsection{Spatial accuracy of scalar transport due to incompressible fluid flow} \label{Nils case}
Finally, we assess the order of accuracy of the advection-diffusion system coupled to an incompressible Navier-Stokes solver.
We consider a circular solid region centered about the point $(\pi, \pi)$ with radius $r = 1.5$.
The penalized momentum, continuity, and advection-diffusion equations are given by
\begin{align}
\D{\rho \u}{t} + \div \rho\u\u &= -\grad p + \div \left[\mu \left(\grad \u + \grad \u^T\right) \right] + \frac{\chi}{\eta}(\ub - \u) + \f, \label{eqn_momentum} \\
\div \u &= 0, \label{eqn_continuity} \\
\D{q}{t} + \left(1 - \chi \right) \left(\u\cdot\grad q\right) &= \div \left[\left\{\kappa \left(1 - \chi\right) + \eta \chi \right\} \grad q \right]+\left(1 - \chi\right) f + \div \left(\chi \vbeta\right) - \chi \div \vbeta. \label{eqn_adv_diff}
\end{align}
Here, $\u(\x,t)$ denotes the fluid velocity, $p(\x,t)$ denotes the fluid pressure, $q(\x,t)$ is a scalar quantity that is passively
transported by the flow, $\f(\x,t)$ denotes the momentum body force term and $\eta$ is the penalization parameter.
The fluid density $\rho$, fluid viscosity $\mu$, and diffusivity $\kappa$ are all set to $1$.
The flux-based VP method is used to impose inhomogeneous Neumann boundary conditions on the surface of the solid for the transported variable $q$ $(\n \cdot \nabla q =  \n \cdot \nabla \qexact)$, while the standard Brinkman penalization method is used to impose Dirichlet boundary conditions for the velocity $(\u = \ub)$. Once again the flux boundary condition for $q$ is spatially varying.

We use the MMS with the following exact steady-state solutions for $\u$, $p$, and $q$:
\begin{align}
&\uexact (\x, t \rightarrow \infty) = \sin(x) \cos (y), \\
& \vexact (\x, t \rightarrow \infty) = - \cos (x) \sin (y), \\
& \pexact(\x, t \rightarrow \infty) = \sin(x)\sin(y), \\
& \qexact (\x, t \rightarrow \infty)= \sin(x)  \sin(y).
\end{align}
These exact solutions are plugged into the \emph{unpenalized} versions of Eq.~\eqref{eqn_momentum} and~\eqref{eqn_adv_diff} in order to determine the required forcing functions $\f$ and $f$.
Note that the imposed boundary condition in the solid region is the steady-state velocity $\u = \ub = \mathbf{u}_\text{exact}$. The fluid and advection-diffusion solvers employed here are second-order accurate in both space and time. All terms in Eqs.~\eqref{eqn_momentum}-\eqref{eqn_adv_diff} are treated implicitly in time, except for the convective terms that are treated explicitly. We refer readers to~\cite{Nangia2019MF} for more details on the spatiotemporal discretization employed in our solvers.
The coupled system is run with a time step size of $\dt = 2 \times 10^{-3}$ (convective CFL is approximately 0.30) until steady-state and error norms are computed
between the exact and numerical solutions within the fluid domain (outside the circular region).
\begin{figure}[]
\centering
%\hspace {25 pt}
\subfigure[Numerical solution of $q$]{
\includegraphics[scale = 0.1]{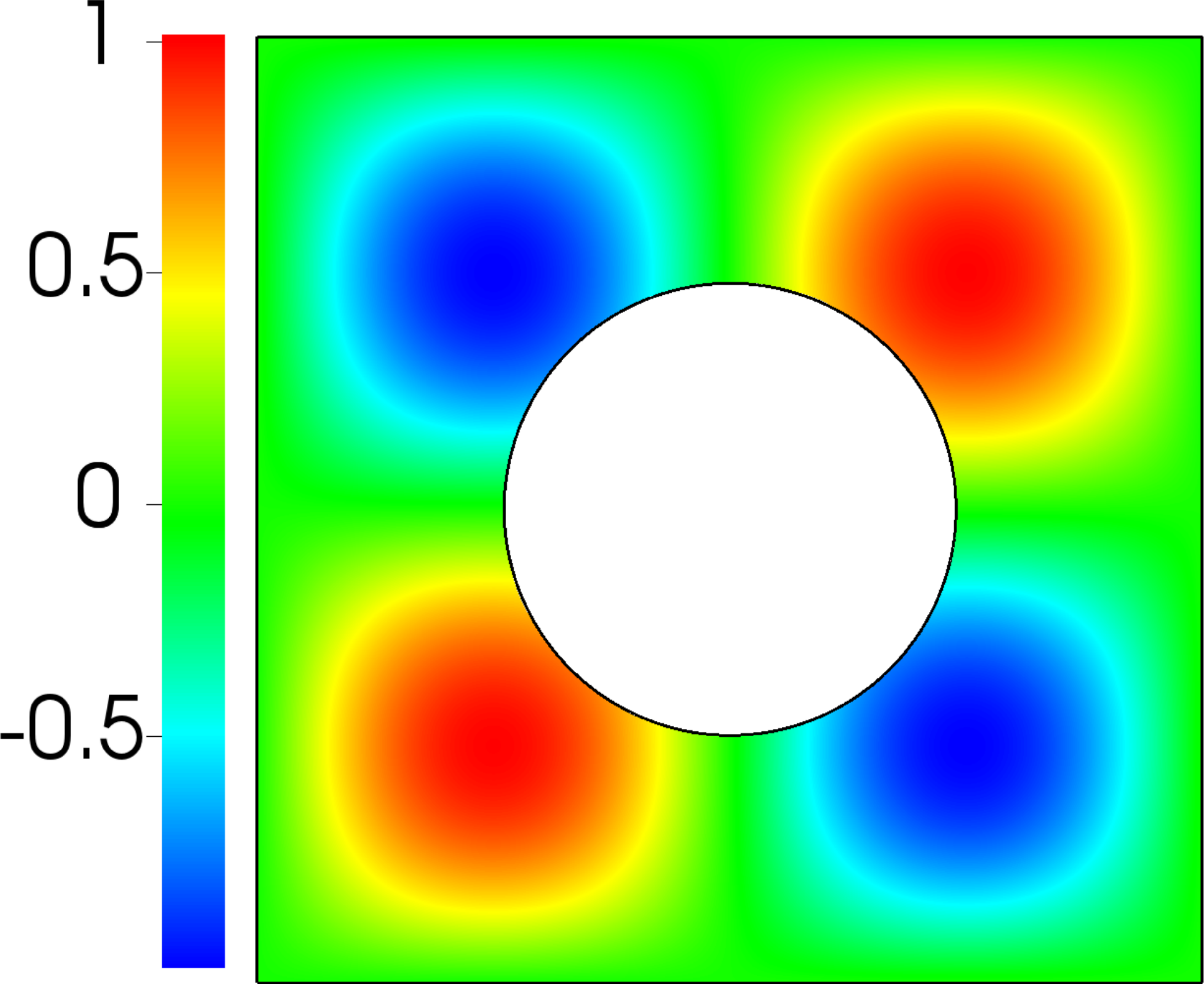}
\label{fig_nils_temperature}
}
\subfigure[Pressure field and velocity vectors]{
\includegraphics[scale = 0.1]{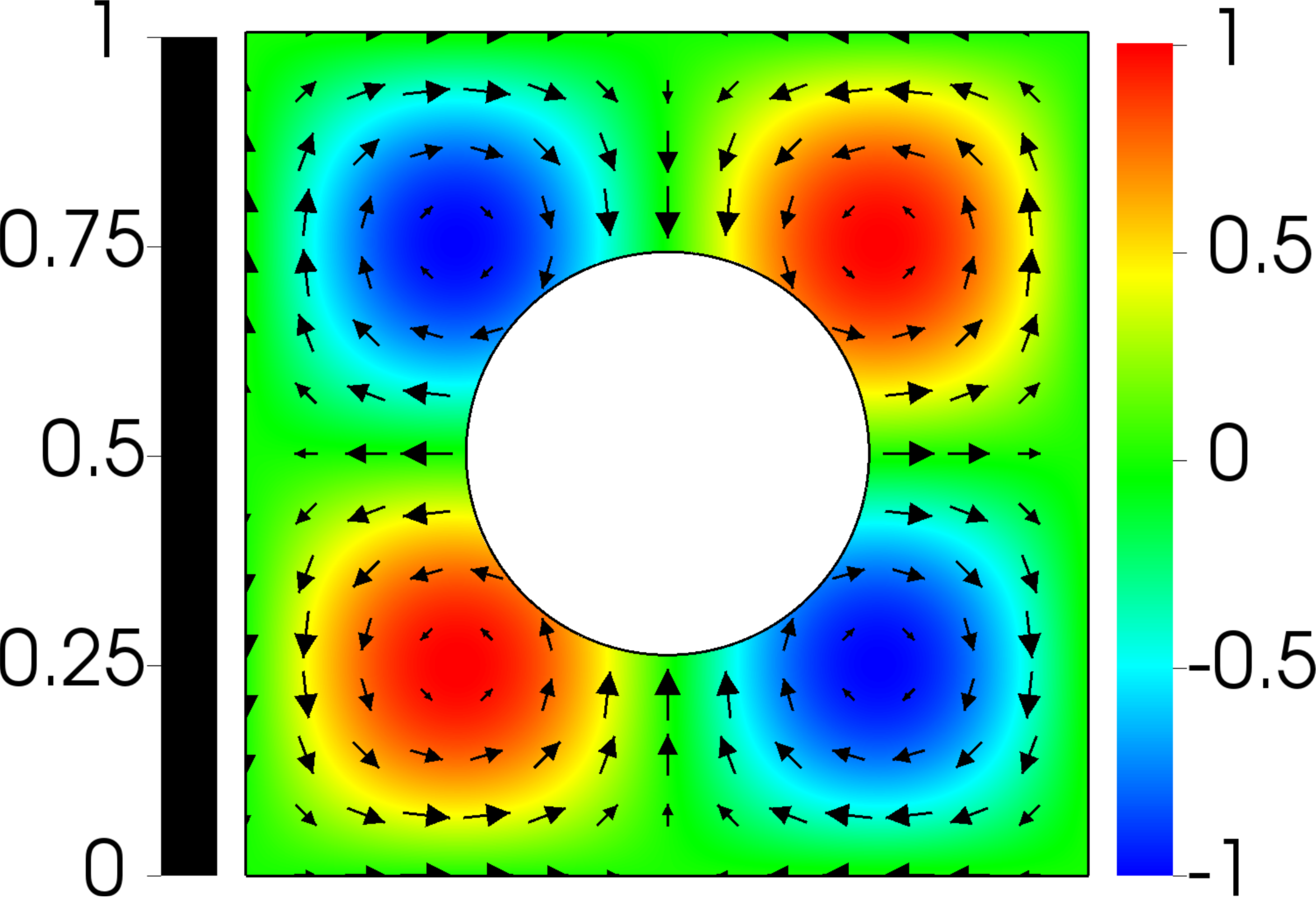}
\label{fig_nils_velocity}
}\\
\subfigure[Spatial convergence rate for $q$]{
\includegraphics[scale = 0.06]{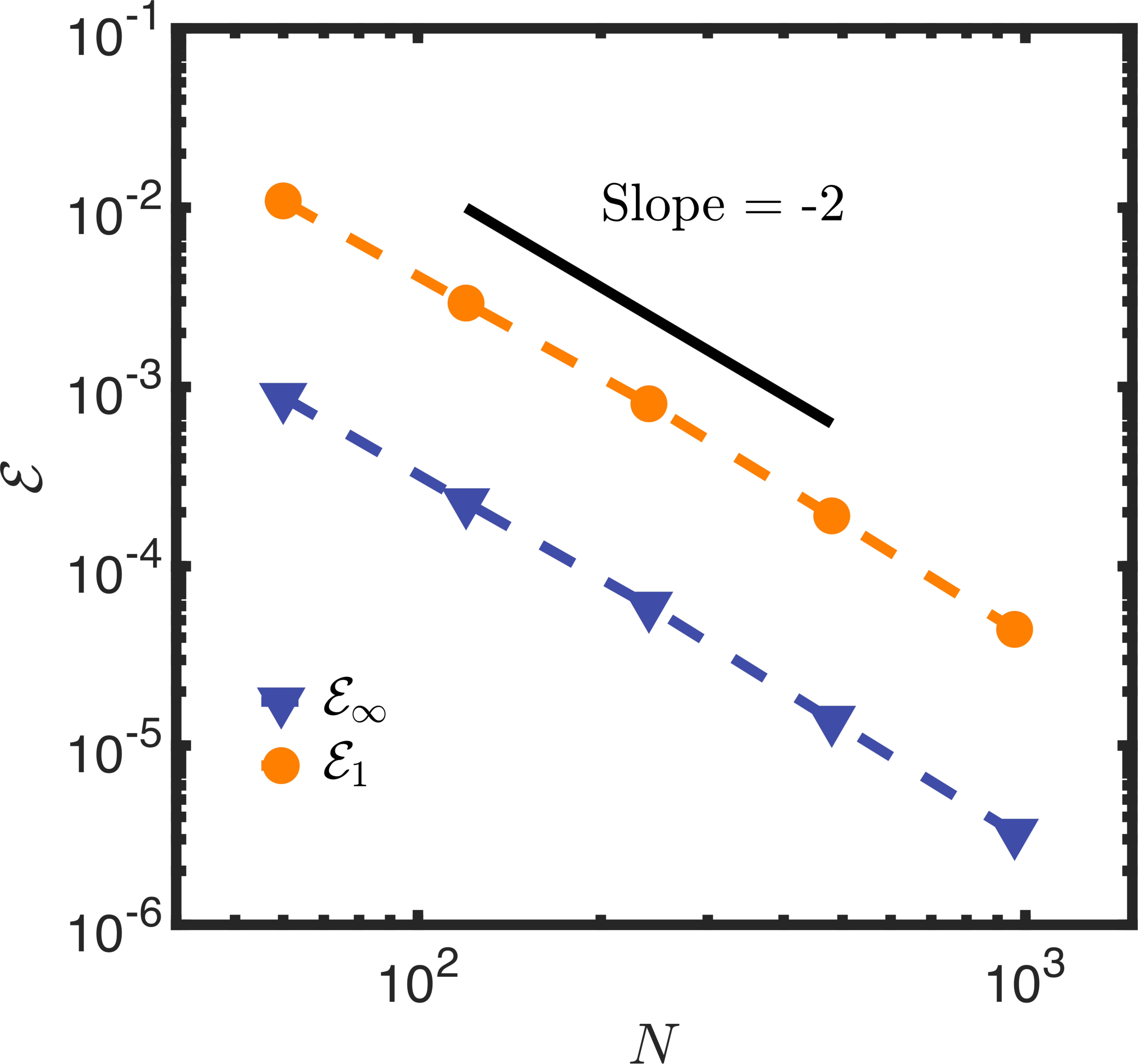}
\label{fig_nils_temperature_convergence}
}
 \subfigure[Spatial convergence rate for $\u$]{
\includegraphics[scale = 0.06]{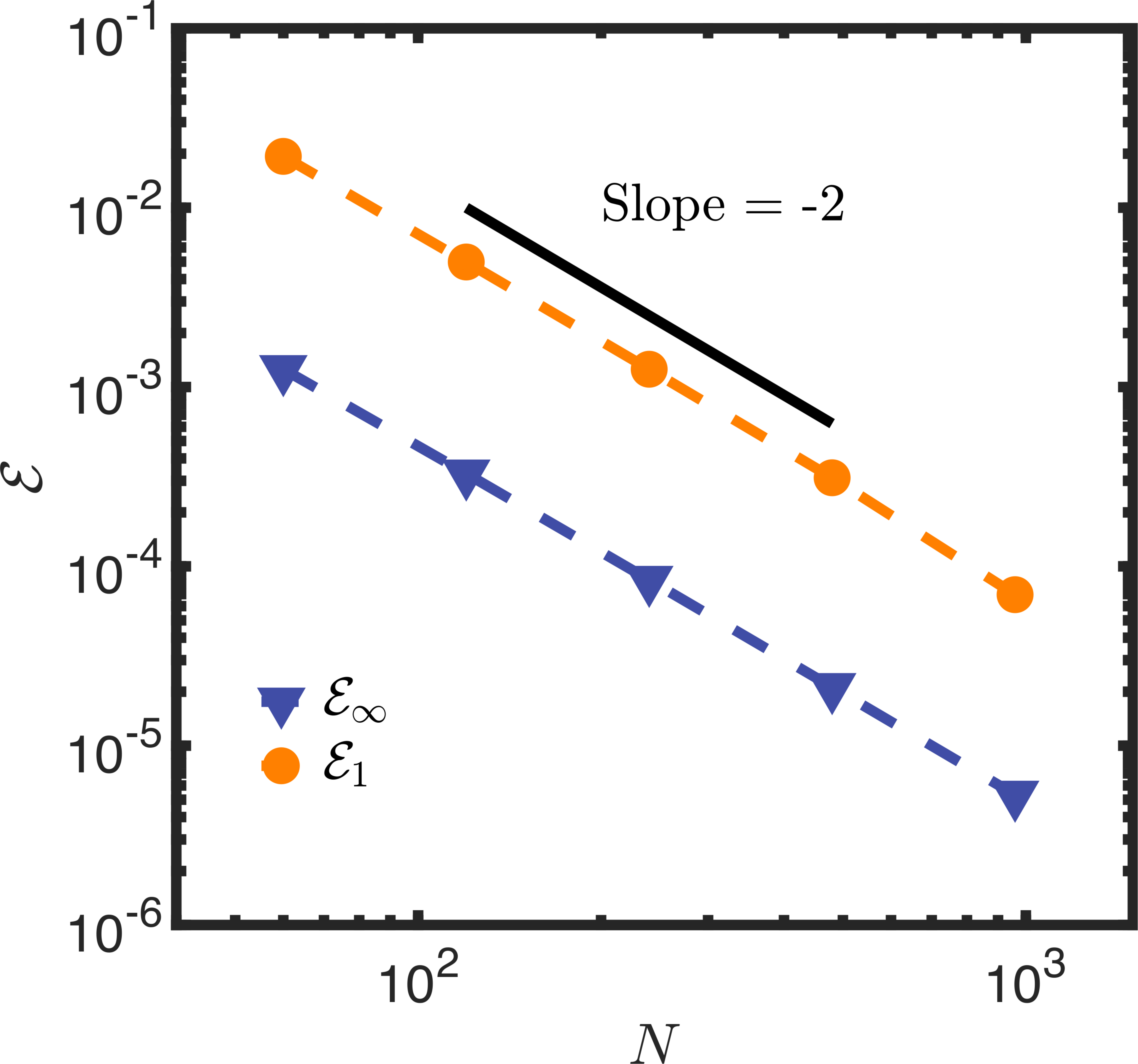}
\label{fig_nils_velocity_convergence}
}
\subfigure[Spatial convergence rate for $p$]{
\includegraphics[scale = 0.06]{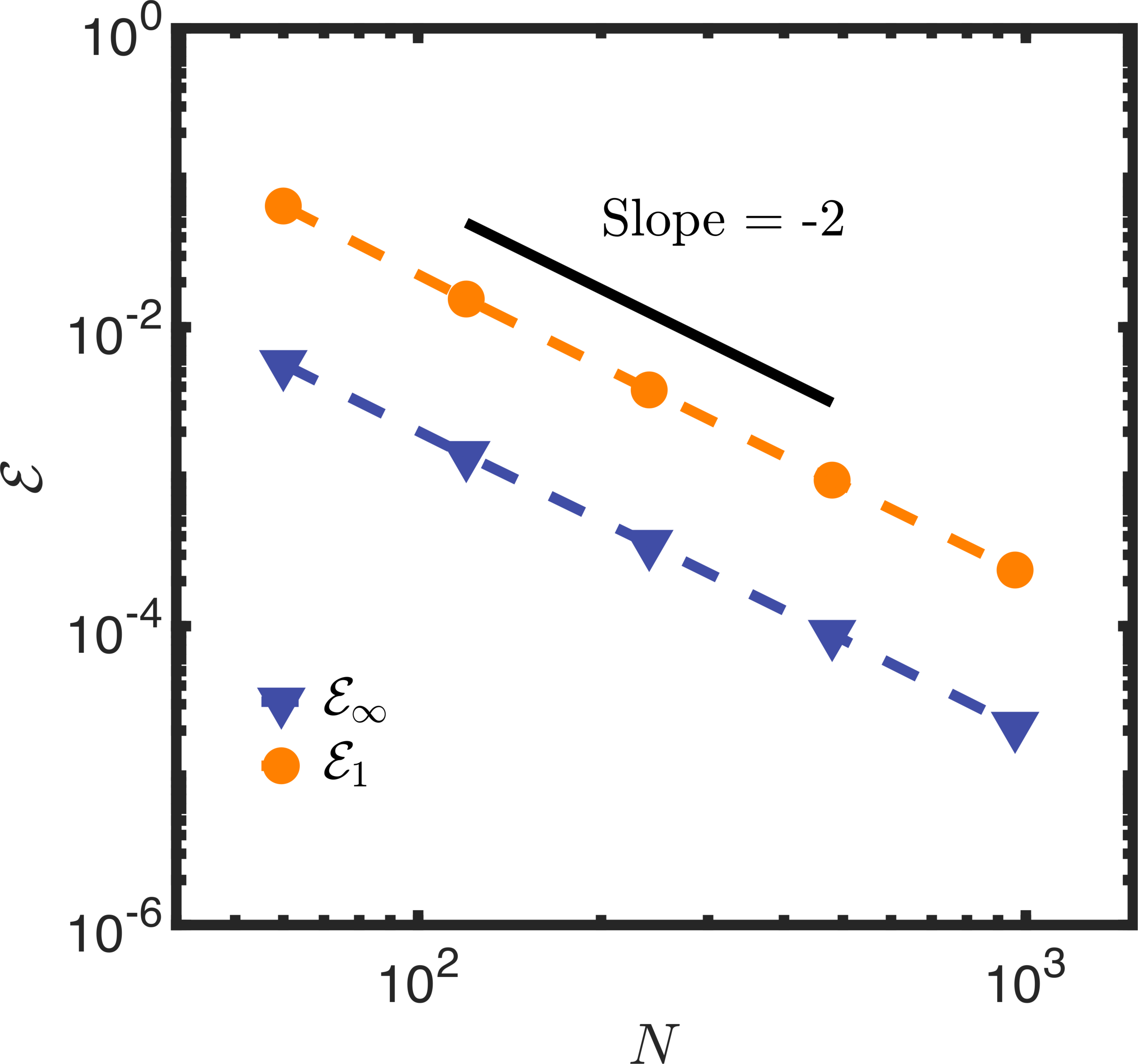}
\label{fig_nils_pressure_convergence}
}
\caption{Advection-diffusion system coupled to a fluid solver problem. Numerical solution of:
\subref{fig_nils_temperature}  Transported quantity $q$;
\subref{fig_nils_velocity} Pressure field and velocity vectors.
Error convergence rate as a function of grid size $N$ for:
\subref{fig_nils_temperature_convergence} $q$,
\subref{fig_nils_velocity_convergence} $\u$, and
\subref{fig_nils_pressure_convergence}$p$. The results shown here are obtained using the continuous indicator function.
}
\label{fig_Nils_case}
\end{figure}

Fig.~\ref{fig_Nils_case} shows the order of convergence for the transported quantity $q$,  velocity $\u$, and pressure $p$. The numerical solutions of $q$ and $p$ are also shown. As observed in the figure, the velocity $\u$,  pressure $p$, and the transported quantity $q$ exhibit second-order convergence rates. This test further corroborates the prior section's results that $\mathcal{O}(2)$ convergence is possible using interfaces that do not conform to the Cartesian grid within the flux-based VP framework.

%% file: Conclusions.tex
%In this letter, we used the method of manufactured solutions to further analyze the spatial order of accuracy of the novel flux-based 
%VP formulation described in~\cite{Sakurai2019}. We demonstrated that the method can be applied to 
%problems involving spatially varying fluxes imposed at the boundary of both simple and relatively complex penalized regions.
%Moreover, we were able to successfully apply the method to problems involving non-periodic boundary conditions on the
%external computational domain. Finally, we investigated the effect of the smoothed masking $\chi$ as it approaches
%a discontinuous indicator function, as well as the influence of different functional forms of the forcing term. Although far from 
%comprehensive, our analysis demonstrates the need to further interrogate the imposition of flux-based
%boundary conditions within the volume penalization framework. Doing so would greatly expand the understanding of 
%Sakurai et al.'s technique and its applicability to more complex multiphysics problems.

In this letter, we used the method of manufactured solution to analyze the spatial order of accuracy of the novel flux-based 
VP formulation described in~\cite{Sakurai2019}. We demonstrated that the flux-based VP method can exhibit second-order spatial convergence even when different flux values are imposed on interfaces that do not conform to the Cartesian grid. We also showed that the convergence rate provided in \cite{Sakurai2019} for some of the cases is not reproducible.  We considered both continuous and discontinuous indicator functions in our test problems.  The two indicator functions yielded similar convergence rate for the problems considered here. Our results  suggest that the flux-based VP approach has a spatial order of accuracy between $\mathcal{O}(1)$ and $\mathcal{O}(2)$, and it depends on the underlying problem/model. The convergence rate cannot simply be deduced \emph{a priori} based on the imposed flux values,  shapes, or grid-conformity of the interfaces, as concluded in Sakurai et al. Further analysis should be carried out to understand the spatial convergence rate of the flux-based VP method. 

We also demonstrated that the method can be applied to problems involving spatially varying flux values on the embedded boundaries. Moreover, cases involving non-periodic boundary conditions on the external computational domain were also considered. Finally, we successfully applied this method to the advection-diffusion equation coupled to an incompressible Navier-Stokes solver, and demonstrated a case in which second-order convergence is achieved for an  (circular) interface that does not conform to the Cartesian grid.

%% file: Supplementary_Material.tex
MATLAB scripts to simulate the 1D Poisson problems of Sec.~\ref{sec_1d_same_flux_bc} and~\ref{sec_1d_different_flux_bc} are included in the supplementary material section. A MATLAB script used to generate signed distance functions for the complex domains considered
in Sec.~\ref{sec_complex_shapes} is also included.  The two dimensional test problems considered in this work can be obtained from the IBAMR Github repository~\cite{IBAMR-web-page}.